\documentclass[aos]{imsart}

\RequirePackage{amsthm,amsmath,amsfonts,amssymb}
\RequirePackage[authoryear]{natbib}
\RequirePackage[colorlinks,citecolor=blue,urlcolor=blue]{hyperref}
\RequirePackage{graphicx}

\startlocaldefs
\numberwithin{equation}{section}
\theoremstyle{plain}
\newtheorem{theorem}{Theorem}
\newtheorem{proposition}{Proposition}
\newtheorem{lemma}{Lemma}
\newtheorem{corollary}{Corollary}
\theoremstyle{remark}
\newtheorem{assumption}{Assumption}

\def\boe{\begin{enumerate}}
\def\eoe{\end{enumerate}}

\newcommand\ca[1]{{\cal{#1}}}
\newcommand\lo[1]{_{\nano{#1}}}
\newcommand\hi[1]{^{\nano{#1}}}

\def\cid{\stackrel{\mbox{\tiny $\cal D$}}\longrightarrow}

\def\proof{\noindent {\sc Proof. }}

\def\cip{\stackrel{\mbox{\tiny $P$}}\rightarrow}

\def\L{{\cal L}}

\def\tsum{\textstyle{\sum}}

\def\trans{^{\mbox{\tiny{\sf  T}}}}

\def\inv{^{\mbox{\tiny $-1$}}}

\newcommand{\indep}{\;\, \rule[0em]{.03em}{.65em} \hspace{-.45em}
\rule[-.05em]{.65em}{.03em} \hspace{-.45em}
\rule[0em]{.03em}{.65em}\;\,}

\def\vec{\mathrm{vec}}

\def\eop{\hfill $\Box$ \\
}
\def\trans{^{\mbox{\tiny{\sf T}}}}

\def\cid{\stackrel{\mbox{\tiny $\cal D$}}\rightarrow}

\def\ali{&\,}

\def\of{{\nano {\circ}}}
\def\nano{\scriptscriptstyle}

\def\real{{\mathbb R}}

\def\tint{\textstyle{\int}}
\def\ka{\kappa}

\def\L2T{L \lo 2 (T)}
\def\L2TX{L \lo 2 (T\lo X)}
\def\L2TX{L \lo 2 (T\lo Y)}
\def\tsum{\textstyle{\sum}}
\def\ali{&\,}

\def\ali{&\,}

\def\eod{\end{document}}

\def\loo#1{_{\nano{\mathrm{\uppercase{#1}}}}}

\newfont{\rsfsten}{rsfs10 scaled 1100}
\newfont{\rsfstena}{rsfs10 scaled 800}
\newfont{\rsfstenb}{rsfs10 scaled 500}

\def\eop{\hfill $\Box$ \\ }
\def\proof{\noindent{\sc{Proof. \ }}}

\def\nano{\scriptscriptstyle}
\def\real{\mathbb R}
\def\vec{\mathrm{vec}}

\def\loo#1{_{\nano{\mathrm{\uppercase{#1}}}}}

\def\nano{\scriptscriptstyle}
\def\inv{\hi{\nano -1}}

\def\nano{\scriptscriptstyle}
\def\of{\mbox{\raisebox{1pt}{$\nano{\circ}$}}}
\def\ka{\kappa}

\def\tsum{{\textstyle{\sum}}}
\def\ali{&\,}

\def\cip{\stackrel{{\nano P}}\rightarrow}
\def\hii#1{\hi{\mbox{\tiny {(\uppercase{#1})}}}}

\def\hiii#1{\hi {(\mathrm{#1})}}
\def\sphere{\mathbb{S}}
\def\tra{\mathrm{tra}}
\newcommand{\RN}[1]{%
  \textup{\uppercase\expandafter{\romannumeral#1}}%
}
\endlocaldefs

\begin{document}

\begin{frontmatter}
\title{A nonparametric test for elliptical distribution based on kernel embedding of probabilities}
\runtitle{Kernel-embedding based test for elliptical distribution}

\begin{aug}
\author[A]{\fnms{Yin}~\snm{Tang}\ead[label=e1]{yqt5219@psu.edu}}
\and
\author[A]{\fnms{Bing}~\snm{Li}\ead[label=e2]{bxl9@psu.edu}}

\address[A]{Department of Statistics, The Pennsylvania State University,\printead[presep={,\ }]{e1,e2}}
\end{aug}

\begin{abstract}
Elliptical distribution is a basic assumption underlying many multivariate statistical methods. For example, in  sufficient dimension reduction and statistical graphical models, this assumption is routinely  imposed to simplify the data dependence structure. Before applying such methods, we need to decide whether the data are elliptically distributed.
Currently existing tests either focus exclusively on spherical distributions, or rely on bootstrap to determine the null distribution, or require specific forms of the alternative distribution. In this paper, we introduce a general nonparametric test for elliptical distribution based on kernel embedding of the probability measure that embodies the two properties that characterize  an elliptical distribution: namely, after centering and rescaling, (1) the direction and length of the random vector are independent, and (2) the directional vector is uniformly distributed on the unit sphere.    We derive the {asymptotic distributions} of the test statistic via von-Mises expansion, develop the sample-level procedure to determine the rejection region, and establish the consistency and validity of the proposed test. {We also develop the concentration bounds of the test statistic, allowing the dimension to grow with the sample size, and further establish the consistency in this high-dimension setting. We compare our method with several existing methods via simulation studies, and} apply our test to  a SENIC dataset with and without a transformation aimed to achieve ellipticity.
\end{abstract}



\end{frontmatter}


\section{Introduction}
\par Elliptical distribution is a widely-used assumption for many statistical and machine learning methods for multivariate data. For example, in sufficient dimension reduction, the elliptical distribution assumption for the predictor is needed for moment-based methods such as  Sliced Inverse Regression
\citep{li1991sliced},  Ordinary Least Squares \citep{li1989regression}, and Iterative Hessian Transformation \citep{Cook2002}. See also \cite{li2018sufficient}. When this assumption is violated, one needs either to perform data transformation or to modify  the inverse-regression methods {as in } \cite{li2009dimension}. Another example is the statistical   graphical model, where a class of methods, such as glasso \citep{yuan2006model} and the transelliptical graphical model \citep{liu2012transelliptical}, require either a Gaussian or an elliptical distribution.  See also \cite{vogel2011elliptical}, which  introduces a class {of} elliptical graphical models as a robust alternative to the Gaussian graphical model.

There are some existing tests on spherical distributions. For example,
\cite{ludwig1991rotationally} proposes a test based on the $L^2$-distance between the empirical distribution  of the data and the distribution function partially specified  by the
spherical assumption. \cite{liang2008necessary} introduces a necessary test by applying the Rosenblatt transformation on each element. The hypothesis in this test is necessary in the sense that it is implied by spherical distribution but does not imply the spherical distribution. \cite{henze2014testing} introduces a test by checking whether the characteristic function is constant over surfaces of spheres centered at the origin. \cite{kariya1977robust} introduces a robust test of a spherical distribution  centered at the origin against an elliptical or a noncentered spherical distribution.  \cite{koltchinskii1998testing} uses the multivariate  distribution and quantile functions to test spherical distributions with unknown centers. All tests above are for the spherical distributions where the covariance matrix is exactly the identity matrix, and most of them focus on the case when we know the center is zero. These methods do not yield direct extensions for testing elliptical distributions where both the mean vector and  {covariance} matrix  are unknown.

There also exist some tests on elliptical distributions. \cite{huffer2007test} gives a test for multivariate normal and elliptical distribution of the data {based on a chi-square statistic after slicing the data}. For the multivariate normal distribution, they derive  the asymptotic null distribution  of their test; for the elliptical distribution, they propose  a bootstrap test without giving a proof of its consistency and validity. \cite{albisetti2020testing} introduces a test based on a Kolmogorov-Smirnov type statistic and uses bootstrap to construct the null distribution. \cite{manzotti2002statistic} proposes a test on whether the standardized directional vector is uniformly distributed on the unit sphere, which is, again, only  a necessary condition for spherical distribution. {\cite{schott2002testing} proposes a Wald-type test based on whether the 4th moments are consistent with an elliptical distribution. So its null hypothesis is not exactly the elliptical distribution. Furthermore, since it requires estimation of the 8th moments, it may not be robust.  \cite{cassart2008optimal} proposes a locally and asymptotically optimal Pseudo-Gaussian test for Fechner-type symmetry, which is wider than the class of elliptical distributions.
}  \cite{babic2019optimal} develops  optimal tests for elliptical distributions against some generalized skew-elliptical alternatives. {Some of the existing methods are summarized in \cite{babic2021rjournal}.}

In this paper, we introduce a nonparametric test for  elliptical distributions based on {Hilbert-space} embedding of a product  probability measure that characterizes an elliptical distribution. The basic idea is the following. It is well known that a random vector $X$ follows a spherical distribution centered at 0 if and only if
\begin{enumerate}
\item its Euclidean norm $\| X \|$ and its direction $X/ \| X \|$ are statistically independent;
\item  the direction vector $X/ \| X \|$ is uniformly distributed on the unit sphere.
\end{enumerate}
{See, for example, \cite{anderson2009introduction,paindaveine2012elliptical}.} This converts testing of sphericity to testing of the two conditions. Since an elliptical distribution can always be transformed into a spherical distribution by a linear map,  we can further develop tests of ellipticity by testing the two conditions for the linearly transformed data. However, since the mean vector and covariance matrix  need to be estimated, we have to take into account the estimation error in this step {when} deriving the asymptotic distribution.

More specifically, {let $U = \|X \|$ and $V = X / \| X \|$. If $X$ has a spherical distribution, then the distribution $P \lo {U,V}$ can be expressed by the product measure $P \lo U \times P \lo V$, where $P \lo V$ is a {\em known} distribution.} We embed this distribution into a reproducing kernel Hilbert space as a cross-covariance operator and compare it against the {kernel embedding} of the fully empirical distribution. The norm of the difference should be small if $X$ has a spherical distribution, and large otherwise. This is the core idea of our method. {One side-note is that we may replace $V$ by its polar coordinate representation, which will significantly simplify the computation.} This procedure has several appealing features: first, the hypothesis is both necessary and sufficient,  that is, $X$ has a spherical distribution if and only if the distance is small; second, it is rather straightforward to go from spherical distribution to elliptical distribution by replacing $X$ with its centered and rescaled version; third, since the test only involves functions of sample moments, its asymptotic null distribution can be relatively easily derived from the infinite-dimensional $\delta$-method, or the von-Mises expansion.

Probability embedding \citep{sriperumbudur2010relation,sriperumbudur2011universality} is a powerful method that has been used in a variety of settings in statistics and machine learning, such as test of independence and the two-sample problem. See, for example, \cite{gretton2005measuring}, \cite{gretton2007kernel}, \cite{gretton2008kernel},  \cite{gretton2009fast}, and \cite{gretton2012kernel}. There is another type  of tests of independence based on the distance covariance; see \cite{szekely2007measuring} and \cite{szekely2009brownian} among others. \cite{sejdinovic2013equivalence} establishes the relation between these two types of tests of independence. Our test of elliptical distribution goes beyond the test of independence between $U$ and $V$, as it must also incorporate the fact that  the distribution of $V$ is known.

{The rest of the paper is organized as follows.  In Section \ref{section:background}, we {lay out the two characterizing properties of an elliptical distribution and construct the probability embedding into a reproducing kernel Hilbert space that embodies the  two characterizing properties.}
{In Sections \ref{section:test statistic} and \ref{section:implement test statistic}, we introduce the test statistic based on the probability embedding  and implement it numerically through coordinate mapping.
In Section \ref{section:null distribution}, for the fixed dimension,  we derive the asymptotic null and alternative distributions of the test statistic via von-Mises expansion, and in Section \ref{section:implement null distribution}, we implement the asymptotic null distribution at the sample level, and establish the validity and consistency of our test.
In Section \ref{section:concentration-bound}, we derive the uniform concentration bounds for our test statistic allowing the dimension to grow with the sample size, and further establish the consistency of our test in this high-dimensional setting.}
In Section \ref{section:simulation}, we conduct simulation studies to demonstrate the usage and effectiveness of the proposed test.
In Section \ref{section:application}, we apply our test to a data example.
{Section \ref{section:discussion} is devoted to some discussions on the choice of kernel functions.

Due to space limit, {some technical lemmas, most of the proofs,  and discussions of some theoretical results} are placed in Appendix A; additional simulation comparisons  are placed in Appendices B and C; further discussions on the choice of kernel functions are placed in Appendix D; and the scatter plot matrices of the dataset in Section \ref{section:application} are placed in Appendix E. All Appendices are  in the online Supplementary Material (\cite{supplementary})}.}
{Our proposed method is implemented in the R package \texttt{KEPTED} (Kernel-Embedding-of-Probability Test for Elliptical Distribution).}

\section{Elliptical distribution and its kernel embedding}\label{section:background}

In this section we introduce the definition of the elliptical distribution and develop two equivalent conditions, one at the level of probability measures and the other at the level of linear operators. The sufficient condition at the operator level is the theoretical basis of our test.

\subsection{Spherical and elliptical distributions}

Let $(\Omega, \ca F, P)$ be a probability space, and $(\Omega \lo X, \ca F \lo X)$ be a measurable space, where $\Omega \lo X$ is a subset of $\real \hi d$, and $\ca F \lo X$ is the Borel $\sigma$-field on $\Omega \lo X$. Let $X: \Omega \to \Omega \lo X$ be a Borel random vector, and let $P \lo X = P \of X \inv$ be the distribution of $X$. {Denote $\lambda$ the Lebesgue measure in $\real \hi d$.} We say that $X$ has a spherical distribution centered at $\mu \in \real \hi d$ if
\begin{enumerate}
\item $P \lo X$ is dominated by the Lebesgue measure in $\real \hi d$ with density $f \lo X = d P \lo X / d \lambda$;
\item $f \lo X(x)$  is a function of $(x-\mu) \trans (x-\mu)$, that is,
\begin{equation}
    f \lo X(x)=h\left((x-\mu) \trans (x-\mu)\right),\label{eq-density-sph}
\end{equation}
for some nonnegative function $h$ satisfying
\begin{equation}
    \int \lo {-\infty} \hi \infty\ldots\int \lo {-\infty} \hi \infty h(x \trans x)dx \lo 1\ldots dx \lo d=1.\label{eq-pdf1}
\end{equation}
\end{enumerate}
Here, $\mu$ is necessarily $EX$ if $X$ is integrable. {Throughout, $\| \cdot \|$ will denote the Euclidean norm (or $L \lo 2$ norm).} Let $U=\|X-\mu\|$ and $V=(X-\mu)/U$. By Theorem 1 of \cite{cambanis1981theory}, $X$ has a spherical distribution if and only if  $U\indep V$ and $V$ has a uniform distribution on the unit sphere $\mathbb{S} \hi {d-1}$ in $\mathbb{R} \hi {d}${, where $\indep$ refers to independent throughout this paper}. See also \cite{eaton1986characterization} and \cite{schmidt2002tail}.

More generally, $X$ is said to have an elliptical distribution centered at $\mu\in\mathbb{R} \hi d$ with a positive definite shape parameter $\Lambda\in\mathbb{R} \hi {d\times d}$ if
\begin{enumerate}
\item $P \lo X$ is dominated by the Lebesgue measure in $\real \hi d$ with density $f \lo X = d P \lo X / d \lambda$;
\item the density of $X$ is a function of $(x-\mu) \trans \Lambda \hi {-1}(x-\mu)$, that is,
\begin{equation*}
    f \lo X(x)=|\Lambda| \hi {-\frac{1}{2}}h\left((x-\mu) \trans \Lambda \hi {-1}(x-\mu)\right),\label{eq-ellip}
\end{equation*}
for some nonnegative function $h$ satisfying (\ref{eq-pdf1}), where $| \Lambda |$ denotes the determinant of $\Lambda$.
\end{enumerate}
A direct corollary is that $Y=\Lambda \hi {-\frac{1}{2}}(X-\mu)$ has a spherical distribution centered at 0. Consequently, if we let
\begin{equation}
    U=\sqrt{(X-\mu) \trans  \Lambda \hi {-1} (X-\mu)},\quad  V=\Lambda \hi {-1/2} (X-\mu)/U,\label{eq-def-uv}
\end{equation}
then  $X$ has an elliptical distribution with center $\mu$ and shape parameter $\Lambda$ if and only if $U\indep V$, and $V$ has a uniform distribution on the unit sphere $\mathbb{S} \hi {d-1}$.

According to Theorem 2.7.2 of \cite{anderson2009introduction}, if the components of $X$ are square-integrable, then
\begin{equation*}
    EX=\mu,\quad \mathrm{var}(X)=\frac{EU \hi 2}{d}\Lambda.
\end{equation*}
Let $\Sigma=\mathrm{var}(X)$. In this paper, we always assume that $X$ has finite mean and variance. {Obviously, neither the dependence between $U$ and $V$ nor the distribution of $V$ will be affected if we replace $\Lambda$ by $\Sigma$ in their definitions in (\ref{eq-def-uv}).} {See \cite{paindaveine2012elliptical} for a detailed discussion.} So, for convenience, we reset $U$ and $V$ to be
\begin{align}\label{eq:u v}
    U=\sqrt{(X-\mu) \trans  \Sigma \hi {-1} (X-\mu)},\quad  V=\Sigma \hi {-1/2} (X-\mu)/U,
\end{align}
for the rest of the paper.
The sufficient and necessary condition for elliptical distribution still applies to the redefined $U$ and $V$, which we record below formally for easy reference.

\begin{proposition} A random vector $X$ has an elliptical distribution if and only if, for $U$ and $V$ defined in (\ref{eq:u v}),
\begin{enumerate}
\item $U \indep V$;
\item $V$ is uniformly distributed in $\mathbb{S} \hi {d-1}$.
\end{enumerate}
\end{proposition}

\subsection{Polar coordinate transformation and equivalent condition}

Since the random vector $V$ only takes values in the unit sphere $\mathbb{S} \hi {d-1}$ in $\real \hi d$, it is more convenient to transform it into a $d-1$ dimensional vector representing the direction of the unit vector $V$ via the polar coordinate system. Specifically, let $v=(v \lo 1, \ldots, v \lo d) \trans \in \mathbb{S} \hi {d-1}$ and let
\begin{align*}
\theta = (\theta \lo 1, \ldots, \theta \lo {d-1}) \trans \in (-\pi/2, \pi/2] \times \cdots \times (-\pi / 2, \pi / 2] \times (-\pi, \pi ] \equiv \Omega \lo \Theta.
\end{align*}

\def\Arctan{\mathrm{Arctan}}
{
For clarity of the subsequent discussion, we first give a definition of the precise meaning of the arc tangent function in the range of $(-\pi, \pi]$. For a real number $r \in \real$, let $\arctan(r)$ be the unique $\theta \in (-\pi/2, \pi/2)$ such that $r = \tan (\theta)$.
Then, for  $(x,y) \in \real \hi 2$, let
\begin{align}\label{eq:def-Arctan}
\Arctan (x,y) =
\begin{cases}
\arctan(y/x), & \mbox{if  $x > 0$}, \\
\arctan(y/x) + \pi, & \mbox{if $x < 0, y \ge 0$ },\\
\arctan(y/x)-\pi, & \mbox{if $x< 0, y < 0$}, \\
\pi/2, & \mbox{if $x = 0$, $y > 0$}, \\
-\pi / 2, & \mbox{if $ x = 0$, $y < 0$}.
\end{cases}
\end{align}
We have written this function as $\Arctan(x,y)$ instead of $\Arctan (y/x)$ because  it is no longer the function of the ratio $y/x$. Some useful properties of the Arctan function \eqref{eq:def-Arctan} are summarized in the Supplementary Materials. The next lemma, whose proof is also placed in the Supplementary Materials, gives the explicit one-to-one correspondence between $v \in \mathbb{S} \hi {d-1}$ and $\theta \in \Omega \lo \Theta$. In the following, we use $S \lo j$ to denote the Euclidean norm of the vector $(v \lo j , \ldots, v \lo d)\trans$.

\begin{lemma}\label{lemma:polar bijection}
The following function from $\Omega \lo \Theta$ to $\sphere \hi {d-1}$
\begin{equation}
    \begin{split}
        v \lo 1&= \sin\theta \lo 1,\\
        v \lo 2&= \cos\theta \lo 1\sin\theta \lo 2,\\
        \vdots&\\
        v \lo {d-1}&= \cos\theta \lo 1\cos\theta \lo 2\ldots\cos\theta \lo {d-2}\sin\theta \lo {d-1},\\
        v \lo d&= \cos\theta \lo 1\cos\theta \lo 2\ldots\cos\theta \lo {d-2}\cos\theta \lo {d-1}.
    \end{split}
    \label{eq-polar-trans}
\end{equation}
is bijective with inverse
\begin{align}
\theta\lo j  =
\begin{cases}
 \Arctan ( S \lo {j+1}, v \lo j )& j = 1, \ldots, d - 2, \label{eq:theta j before} \\
\Arctan(v \lo d, v \lo {d-1}) & j = d - 1.
\end{cases}
\end{align}
\end{lemma}
We will denote the first function \eqref{eq-polar-trans} as $v=\rho (\theta)$ and the second function \eqref{eq:theta j before} as $\theta = g(v)$.
}
 Then
\begin{align*}
x = u v = u  \rho  (\theta) \equiv  \tau  (u, \theta ),
\end{align*}
Evidently $ \tau $ is an invertible function, and the joint distribution of $U, \Theta$ can be written as
\begin{align}
f \lo {U \Theta}(u, \theta) = f \lo X ( \tau  (u, \theta)) \left| \tfrac{\partial  \tau   (u, \theta)} {\partial (u, \theta \trans)} \right|. \label{eq: utheta dist}
\end{align}
According to \cite{anderson2009introduction}, the Jacobian on the right-hand side above is
\begin{align}
 \left| \tfrac{\partial  \tau   (u, \theta)} {\partial (u, \theta \trans)} \right| = u \hi {d-1}\cos \hi {d-2}\theta \lo 1\cos \hi {d-3}\theta \lo 2\ldots\cos\theta \lo {d-2}.
\label{eq: utheta jacobi}
\end{align}
Hence, {by (\ref{eq-density-sph}), (\ref{eq: utheta dist}) and (\ref{eq: utheta jacobi})}, the joint p.d.f. of $(U, \Theta)$ is
\begin{align}
\begin{split}
    f \lo {U, \Theta} (u, \theta) = \ali u \hi {d-1}\cos \hi {d-2}\theta \lo 1\cos \hi {d-3}\theta \lo 2\ldots\cos\theta \lo {d-2}h(u \hi 2), \quad (u, \theta) \in \Omega \lo U \times \Omega \lo \Theta,
\end{split}\label{eq-density-sph-polar}
\end{align}
where $\Omega \lo U= [0,\infty)$. This relation implies (i) $U$ and $\Theta$ are independent, and (ii) $\Theta$ has a known distribution. We summarize this result as the following proposition.

\begin{proposition}\label{proposition:n s condition 1} If $(U, \Theta)= \tau  \inv (X)$, then $X$ has a spherical distribution {centered at 0} if and only if
\begin{enumerate}
\item $U \indep \Theta$, or equivalently $P \lo {U, \Theta} = P \lo U \times P \lo \Theta$;
\item $P \lo \Theta$ has p.d.f.
$
f \lo \Theta (\theta) = c  \cos \hi {d-2} (\theta \lo 1) \cos \hi {d-3}( \theta \lo 2) \cdots \cos( \theta \lo {d-2}),
$
where  $\theta \in \Omega \lo \Theta$ and
\begin{align*}
c  = \left(\tint \lo {\Omega \lo \Theta} \, \cos \hi {d-2} (\theta \lo 1) \cos \hi {d-3}( \theta \lo 2) \cdots \cos( \theta \lo {d-2}) \, d \theta \right)\inv .
\end{align*}
\end{enumerate}
\end{proposition}

\subsection{Kernel embedding of $P \lo U \times P \lo \Theta$}\label{section:kernel-embedding}

Let $\ka \lo U: \Omega \lo U \times \Omega \lo U \to \real$ and $\ka \lo \Theta: \Omega \lo \Theta \times \Omega \lo \Theta \to \real$ be positive definite kernels, and let $\ca H \lo U$ and $\ca H \lo \Theta$ be the reproducing kernel Hilbert space (RKHS) generated by $\ka \lo U$ and $\ka \lo \Theta$. Let $\ca B$ be the collection of linear operators
$\{ f \otimes g: f \in \ca H \lo U, g \in \ca H \lo \Theta \}$. Thus, each member of $\ca B$ is a linear operator mapping from $\ca H \lo \Theta$ to $\ca H \lo U$ such that, for any $h \in \ca H \lo \Theta$, $(f \otimes g ) (h) = f \langle g, h \rangle \lo {\ca H \lo \Theta}$. Let $\ca S$ be the linear span of $\ca B$, consisting of finite linear combinations of members of $\ca B$ with real coefficients. Define in $\ca S$  the inner product
\begin{align*}
\ali \langle \alpha \lo 1 (f \lo 1 \otimes g \lo 1) + \cdots + \alpha \lo r (f \lo r \otimes g \lo r),
 \tilde \alpha \lo 1 ( \tilde f \lo 1 \otimes  \tilde g \lo 1) + \cdots +  \tilde \alpha \lo r ( \tilde f \lo s \otimes  \tilde g \lo s) \rangle \\
\ali \hspace{1in} =
\tsum \lo {i=1} \hi r\tsum \lo {j=1} \hi s \alpha \lo i \tilde \alpha \lo j \langle f \lo i, \tilde f \lo i \rangle \lo {\ca H \lo U} \, \langle g \lo j, \tilde g \lo j \rangle \lo {\ca H \lo \Theta}.
\end{align*}
Endowed with this inner product, $\ca S$ is an inner product space; its completion as a Hilbert space is the tensor product space $\ca H \lo U \otimes \ca H \lo \Theta$.

Let $\ca F \lo U$ and $\ca F \lo \Theta$ be the Borel $\sigma$-fields on $\Omega \lo U$ and $\Omega \lo \Theta$, and let $\ca F \lo U \times \ca F \lo \Theta$ be the product $\sigma$-field. Abbreviate $\Omega \lo U \times \Omega \lo \Theta$ and $\ca F \lo U \times \ca F \lo \Theta$ by $\Omega \lo {U, \Theta}$ and $\ca F \lo {U, \Theta}$. Let $\ca M ( \Omega \lo {U, \Theta}, \ca F \lo {U, \Theta})$ denote the class of all probability measures on $( \Omega \lo {U, \Theta}, \ca F \lo {U, \Theta} )$. We want to find an injective mapping from $\ca M ( \Omega \lo {U, \Theta}, \ca F \lo {U, \Theta})$ to $\ca H \lo U \otimes \ca H \lo \Theta$ so that testing equality of two measures in $\ca M ( \Omega \lo {U, \Theta}, \ca F \lo {U, \Theta} )$ is equivalent to testing the equality of two operators in $\ca H \lo U \otimes \ca H \lo \Theta$. {Such a mapping is provided by the following theorem. {Recall that a kernel $\kappa$ is characteristic if the mapping $P \mapsto \int \kappa (\cdot, X) d P$ is injective. Here, the dot notation $\int \kappa (\cdot, X) d P$ simply means the function $x \mapsto  \int \kappa (x, X) d P$. So $\kappa$ being characteristic means if $\int \kappa (x, X) d P \lo 1 = \int \kappa (x, X) d P \lo 2$ for all $x$ then $P \lo 1= P \lo 2$.}

\begin{theorem}\label{thm:joint characteristic}
 If  $\ka \lo U$ and $\ka \lo \Theta$ are characteristic, then so is tensor product kernel   $\ka \lo U\otimes\ka \lo \Theta$; that is, the mapping
\begin{align*}
\ca M ( \Omega \lo {U, \Theta}, \ca F \lo {U, \Theta} ) \to \ca H \lo U \otimes \ca H \lo \Theta, \quad P \lo {U, \Theta} \mapsto \tint \lo {\Omega \lo {U, \Theta}} \ka \lo U (\cdot, U) \otimes \ka \lo \Theta (\cdot, \Theta) d P \lo {U, \Theta}
\end{align*}
is injective.
\end{theorem}

{The dot notation $\ka \lo U (\cdot, U)$ means the function $u \mapsto \ka \lo U (u, U)$; the same applies to $\ka \lo \Theta (\cdot, \Theta)$. So the operator $\ka \lo U (\cdot, U) \otimes \ka \lo \Theta (\cdot, \Theta)$ maps a function $f \in \ca H \lo \Theta$ to $f(\Theta) \ka \lo U (\cdot, U)$, which is a member of $\ca H \lo U$. Correspondingly, the operator $\int \lo {\Omega \lo {U, \Theta}} \ka \lo U (\cdot, U) \otimes \ka \lo \Theta (\cdot, \Theta) d P \lo {U, \Theta}$ maps a function $f \in \ca H \lo \Theta$ to
$\int \lo {\Omega \lo {U, \Theta}} \ka \lo U (\cdot, U) f(\Theta) d P \lo {U, \Theta}$, which is a member of $ \ca H \lo U$. This type of dot notations will be used throughout the rest of the paper. }
The proof of Theorem \ref{thm:joint characteristic} can be found in Theorem 4 of \cite{szabo2018characteristic}. As a result, we only need to guarantee that both $\ka \lo U$ and $\ka \lo \Theta$ are characteristic kernels.
In fact, the Gaussian radial basis kernels, which in our context are
\begin{equation}\label{eq:gaussian kernel}
\kappa \lo U(u \lo 1,u \lo 2)=\exp[-\gamma \lo U(u \lo 1-u \lo 2) \hi 2],\quad \kappa \lo \Theta(\theta \lo 1,\theta \lo 2)=\exp(-\gamma \lo \Theta\|\theta \lo 1-\theta \lo 2\| \hi 2),
\end{equation}
for some $\gamma \lo U,\gamma \lo \Theta > 0$, indeed satisfy the condition in Theorem \ref{thm:joint characteristic}: it is shown in Theorem 3.2 of \cite{guella2021gaussian} that a Gaussian radial basis kernel is integrally strictly positive definite, and it is further shown that an integrally strictly positive definite kernel is characteristic (see, for example, \cite{sriperumbudur2011universality, fukumizu2009kernel, sriperumbudur2010hilbert}).
}

Theorem \ref{thm:joint characteristic} implies the following equivalence which is the basis of our test of ellipticity.

\begin{corollary}\label{corollary:second equivalence} Suppose
\begin{enumerate}
\item $X$ is a random vector in $\real \hi d$ with mean $\mu$ and covariance matrix $\Sigma$;
\item $(U, \Theta) = \tau  \inv ( \Sigma \hi {-1/2} (X - \mu ))$;
\item $\ka \lo U$ and $\ka \lo \Theta$ are {characteristic} kernels.
\end{enumerate}
Then $X$ has an elliptical distribution with parameters $\mu$ and $\Sigma$ if and only if
\begin{align*}
\int \lo {\Omega \lo {U, \Theta}}[\ka \lo U (\cdot, U) \otimes \ka \lo \Theta (\cdot, \Theta) ] d P \lo {U, \Theta} =
\int \lo {\Omega \lo {U}}\ka \lo U (\cdot, U)  d P \lo {U}
 \otimes \int \lo {\Omega \lo \Theta} \ka \lo \Theta (\cdot, \Theta) d { P \lo {0}}
\end{align*}
where {$P \lo 0$ is the known true $P \lo {\Theta}$} with its form  given in Proposition \ref{proposition:n s condition 1}.
\end{corollary}

\section{Construction of Test Statistic}\label{section:test statistic}

In this section we construct our test statistic based on Corollary \ref{corollary:second equivalence}.
Let $X \lo 1, \ldots, X \lo n $ be an i.i.d. sample of $X$. For a function $A(x)$, we use $E \lo n A(X)$ to denote the sample average of $A(X \lo 1), \ldots, A (X \lo n)$. Let
$\hat \mu$ and $\hat \Sigma$ denote the sample mean and sample variance; that is, $\hat \mu = E \lo n (X)$ and $\hat \Sigma = E \lo n [(X- \hat \mu ) (X - \hat \mu )\trans]$.  For any $x \in \Omega \lo X$, let
{
\begin{align*}
W (x)  = \ali   \Sigma \hi {-1/2} (  x  - \mu ) , \hspace{.25in} \hat W  (x) =  \hat \Sigma \hi {-1/2} (  x  - \hat \mu ), \\
U (x)  = \ali  \| W(x) \|, \hspace{.54in} \hat U (x)  =  \| \hat W (x) \|, \\
V  (x) =  \ali W(x) / U (x) , \hspace{.28in} \hat V (x)  =  \hat W (x) /\hat  U (x), \\
\Theta(x) = \ali g(V(x)),\hspace{.54in} \hat{\Theta}(x)=g(\hat{V}(x)),
\end{align*}
where $g(\cdot)$ is the polar coordinate transformation as defined in \eqref{eq-polar-trans}. When no ambiguity is likely, we abbreviate $W(X)$ by $W$,  $W(X \lo i)$ by $W \lo i$, and $\hat W (X \lo i)$ by   $\hat W \lo i$. The same applies to $U$, $V$ and $\Theta$.}
Let
{
\begin{equation}
\Sigma \lo {U \Theta}=E[\kappa \lo U(\cdot,U )\otimes\kappa \lo \Theta(\cdot,\Theta )]-E[\kappa \lo U(\cdot,U )]\otimes E \lo 0 [\kappa \lo \Theta(\cdot,\Theta)]\label{eq-pop-op-def}
\end{equation}
where $E \lo 0$ refers to the expectation with respect to $P \lo 0$.
Note that $\Sigma \lo {U\Theta}$ is different from a usual cross-covariance operator between two random variables, because $P \lo \Theta$ is given the specific form $P \lo 0$.
By Corollary \ref{corollary:second equivalence}, $X$ has an elliptical distribution if and only if $\Sigma \lo {U \Theta}=0$, which implies $U \indep \Theta$ and $P \lo \Theta$ has the form given in Proposition \ref{proposition:n s condition 1}. Thus, our goal is to test the hypothesis
\begin{align}\label{eq:targeted hypothesis}
H \lo 0: \Sigma \lo {U \Theta} = 0.
\end{align} }
Note that this is not merely a test of independence between $U$ and $\Theta$, because $P \lo \Theta$ has a known, specific form.

Let
\begin{align}
\breve{\Sigma} \lo {U \Theta} = E \lo n [ \kappa \lo U(\cdot,\hat U )\otimes\kappa \lo \Theta(\cdot,\hat \Theta )] - E \lo n [ \kappa \lo U(\cdot,\hat U )] \otimes{ E \lo 0[ \kappa \lo \Theta(\cdot, \Theta )]},\label{eq-sample-op-def}
\end{align}
where, for example, the first term on the right is the sample average of
\begin{align*}
\{ \kappa \lo U(\cdot,\hat U (X \lo i)  )\otimes\kappa \lo \Theta(\cdot,\hat \Theta (X \lo i)): i = 1, \ldots, n \}.
\end{align*}
Since the last term is an integral with respect to the known distribution $P \lo 0$, this operator is {\em not} the usual sample estimate of the {cross-covariance operator} $\Sigma \lo {U \Theta}$. That is why we denote it by $\breve{\Sigma} \lo {U \Theta} $ instead of $\hat {\Sigma} \lo {U \Theta} $.
{The operator $\breve \Sigma \lo {U \Theta}$ is a mapping from $\ca H \lo \Theta$ to $\ca H \lo U$: for a given  $f \in \ca H \lo \Theta$, $\breve{\Sigma} \lo {U \Theta} f$ is the function
\begin{align*}
E \lo n [\ka \lo U (\cdot, \hat U) f(\hat \Theta)] - E \lo n [ \ka \lo U (\cdot, \hat U) ] E \lo 0 [f(\Theta)],
\end{align*}
which is a member of $\ca H \lo U$.}

For convenience, let $\tilde{\kappa} \lo \Theta(\cdot,\theta)=\kappa \lo \Theta(\cdot,\theta)- \int \kappa \lo \Theta (\cdot,\theta)d {P \lo 0} (\theta) $ denote the centered  kernel function in $\mathcal{H}\lo \Theta$. Then, \eqref{eq-pop-op-def} and \eqref{eq-sample-op-def} can be simplified as
\begin{align*}
\Sigma \lo {U \Theta}=E[\kappa \lo U(\cdot,U)\otimes\tilde{\kappa} \lo \Theta(\cdot,\Theta)],\quad \breve{\Sigma} \lo {U \Theta} = E \lo n [ \kappa \lo U(\cdot,\hat U )\otimes\tilde{\kappa} \lo \Theta(\cdot,\hat \Theta )].
\end{align*}
We use
\begin{align*}
T \lo n =n\|\breve{\Sigma} \lo {U \Theta}\| \loo {HS} \hi 2
\end{align*}
as our test statistic for the hypothesis (\ref{eq:targeted hypothesis}).
{The squared Hilbert-Schmidt norm  $\| \breve \Sigma \lo {U\Theta} \| \loo {HS} \hi 2$ is the sum of eigenvalues of the operator $\breve \Sigma \lo {U\Theta} \breve \Sigma \lo {U\Theta} \hi *$, where $\breve \Sigma \lo {U\Theta} \hi *$ is the adjoint operator of $\breve \Sigma \lo {U\Theta}$. This is analogous to the Frobenius norm of a matrix. }

By our definition, $T \lo n / n $ is close to $ \|\Sigma \lo {U \Theta}\| \loo {HS} \hi 2$. {An equivalent definition of $\|\Sigma \lo {U \Theta}\| \loo {HS}$ is
\begin{align*}
\|\Sigma \lo {U \Theta}\| \loo {HS}
&=& \|E[\kappa \lo U(\cdot,U )\otimes\kappa \lo \Theta(\cdot,\Theta )]-E[\kappa \lo U(\cdot,U )]\otimes E \lo 0 [\kappa \lo \Theta(\cdot,\Theta)]\| \loo {HS}  \\
&=& \sup \lo {f \lo 1\otimes f \lo 2 \in \mathcal{F}} | E \lo {P \lo {U \Theta}}[ f \lo 1(U) f \lo 2(\Theta)]-E \lo {P \lo U} [f \lo 1(U)] E \lo {P \lo 0} [f \lo 2(\Theta)]|
\end{align*}
where $\mathcal{F} = \{ f \lo 1 \otimes f \lo 2 \in \mathcal{H} \lo U \otimes \mathcal{H} \lo \Theta: \| f \lo 1 \otimes f \lo 2 \| \le 1 \}$, and $P \lo 0$ is the true distribution of $\Theta$ as defined in Proposition  \ref{proposition:n s condition 1}. In this way, we can interpret $\|\Sigma \lo {U \Theta}\| \loo {HS} $ as a ``distance'' between $P \lo {U \Theta}$ and the ``closest'' elliptically symmetric distribution $P \lo U \times P \lo 0$. Here, ``closest'' can be interpreted as keeping the marginal distribution $P \lo U$ unchanged and using it to construct an elliptical distribution.
 Therefore,} intuitively, $T \lo n / n$ is small if $X$ has an elliptical distribution. On the other hand, if $X$ does not have an elliptical distribution, then  either $U$ and $\Theta$ are not independent, or the marginal distribution of $\Theta$ is not the one given by Proposition \ref{proposition:n s condition 1}. In either case $T \lo n / n$ will not be small.

{
As we will see in the next section, $T \lo n /n $ can be re-expressed as
\begin{align*}
n \hi {-2} \tsum \lo {i,j} \kappa \lo U (\hat U \lo i, \hat U \lo j) \langle \tilde \kappa \lo \Theta (\cdot, \hat \Theta \lo i) ,  \tilde \kappa \lo \Theta (\cdot, \hat \Theta \lo j) \rangle \lo {\ca H \lo \Theta}.
\end{align*}
If $U$ and $\Theta$ are independent, then the above is approximately
\begin{align*}
\left(n \hi {-2} \tsum \lo {i,j} \kappa \lo U (\hat U \lo i, \hat U \lo j)\right) \left(
n \hi {-2} \tsum \lo {i,j} \langle \tilde \kappa \lo \Theta (\cdot, \hat \Theta \lo i) ,  \tilde \kappa \lo \Theta (\cdot, \hat \Theta \lo j) \rangle \lo {\ca H \lo \Theta}\right).
\end{align*}
If, furthermore, $\Theta$ is distributed as $P \lo 0$, then the second term will be near 0. Hence $T \lo n/n$ will be near 0 if and only if both $\Theta \indep U$ and $\Theta \sim P \lo 0$ hold; that is, $X$ has an elliptical distribution.
}

\section{Computing the test statistic}\label{section:implement test statistic}
\subsection{Coordinate mapping}\label{subsection:coordinate mapping}
The implementation of the test at the sample level relies on coordinate representation of linear operators.  Let $\ca H$ be an $r$-dimensional space with a basis $\ca B = \{ b \lo 1, \ldots, b \lo r\}$. Then every member $f$ of $\ca H$ can be represented as $c \lo 1 b \lo 1 + \cdots + c \lo r b \lo r$. The mapping from $\ca H$ to $\real \hi r$ defined by {$C (f) = (c \lo 1, \ldots, c \lo r) \trans$} is called the coordinate mapping. Let $G$ be the Gram matrix $\{ \langle b \lo i, b \lo j \rangle \lo {\ca H}\} \lo {i,j=1} \hi r$. Then, for any $f, g \in \ca H$,
\begin{align*}
\langle f, g \rangle \lo {\ca H } = \tsum \lo {i=1} \hi r  \tsum \lo {j=1} \hi r C(f) \lo i C(g) \lo j \langle b \lo i, b \lo j \rangle \lo {\ca H} = C(f)\trans G C (g).
\end{align*}
In other words, if we let $\real \hi r (G)$ represent the Hilbert space consisting of the vector space $\real \hi r$ along with the inner product {$\langle a, b \rangle = a\trans G b$}, then  $C: \ca H \to \real \hi r (G)$ is an isomorphism. Let $A: \ca H \to \ca H$ be a self-adjoint operator. An eigenvalue of $A$ is defined by the following relations
\begin{align*}
A f = \lambda f, \quad \langle f, f \rangle \lo {\ca H} = 1,
\end{align*}
or equivalently,
\begin{align*}
C A C \hi * C  f = \lambda C f, \quad (C f) \trans G C f = 1,
\end{align*}
where, $C \hi *$ stands for the  adjoint operator of $C$. Note, again, that $*$ and $\star$ denote different concepts.
Letting $v = G \hi {1/2} Cf$, the above can be restated as
\begin{align*}
G \hi {1/2} C A C \hi * G \hi {-1/2}   v = \lambda v , \quad    v \trans v = 1.
\end{align*}
In other words, a number is an eigenvalue of the operator $A$ if and only if it is an eigenvalue of the matrix $G \hi {1/2} C A C \hi * G \hi {-1/2}$, which can be shown to be a symmetric matrix. In particular,
\begin{align*}
\tra (A) = \tra ( G \hi {1/2} C A C \hi * G \hi {-1/2} )= \tra ( C A C \hi *  ),
\end{align*}
where $\tra$ on the left represents the trace of a linear operator, in the middle and on the right represent  that of a matrix.
This identity allows us to express the Hilbert-Schmidt norm of an operator as the trace of a matrix.

Let $e \lo i$ represent the $i$th column of the identity matrix $I \lo r$. Since $e \lo i = C b \lo i$,  $C A C \hi *$ is simply the $n \times n$ matrix
\begin{align}\label{eq:CAC}
C A C \hi * (e \lo 1, \ldots, e \lo p) = C A C \hi * (C b \lo 1, \ldots, C b \lo p) =( C  A b \lo 1  , \ldots, C A b \lo r  ).
\end{align}

\subsection{Test statistic}
At the sample level, $\ca H \lo U$ and $\ca H \lo \Theta$ are spaces spanned by the two bases
\begin{equation*}
\mathcal{B} \lo U=\{\kappa \lo U(\cdot,\hat{U} \lo i):i=1,\ldots,n\},\quad \mathcal{B} \lo \Theta=\{\kappa \lo \Theta(\cdot,\hat{\Theta} \lo i):i=1,\ldots,n\}
\end{equation*}
respectively. Let $K \lo U$ and $\tilde K \lo \Theta$ be the $n \times n$ matrices whose $(i,j)$th entries are
\begin{equation}\label{eq:K and K}
({K} \lo U) \lo {ij}=  \kappa \lo U(\hat{U} \lo i,\hat{U} \lo j), \quad
(\tilde{{K}} \lo \Theta) \lo {ij}= \langle \tilde \ka \lo \Theta ( \cdot, \hat \Theta \lo i), \tilde \ka \lo \Theta ( \cdot, \hat \Theta \lo j) \rangle \lo {\ca H \lo \Theta}
\end{equation}
respectively. Note that $K \lo U$ is simply the Gram matrix of $\ca B \lo U$, and  $(\tilde{{K}} \lo \Theta) \lo {ij}$ can be expanded as
\begin{align*}
 {\kappa} \lo \Theta(\hat{\Theta} \lo i, \hat \Theta \lo j   ) - \tint  \ka \lo \Theta (\hat \Theta \lo i, \theta ) d P \lo  0 (\theta) - \tint  \ka \lo \Theta (\hat \Theta \lo j, \theta ) d P \lo 0 (\theta) + \tint\tint \ka \lo \Theta (\theta, \theta') d (P \lo 0 \times P \lo 0) (\theta, \theta').
\end{align*}
Let $C: \ca H \lo U \to \real \hi n$ be the coordinate mapping. Our goal is to compute
\begin{align*}
T \lo n = n \| \breve \Sigma \lo {U \Theta } \| \loo{HS} \hi 2 = n \, \tra (  \breve \Sigma \lo {U \Theta }  \breve \Sigma \lo {U \Theta } \hi * ).
\end{align*}
By the discussion in Section \ref{subsection:coordinate mapping}, we have
\begin{align*}
\tra (  \breve \Sigma \lo {U \Theta }  \breve \Sigma \lo {U \Theta } \hi * )
=\tra ( C \breve \Sigma \lo {U \Theta }  \breve \Sigma \lo {U \Theta } \hi * C \hi * ).
\end{align*}
The next proposition gives the coordinate of  $\breve \Sigma \lo {U \Theta }  \breve \Sigma \lo {U \Theta } \hi *$.

\begin{proposition}\label{prop: coordinate mapping}
 If  $K \lo U$ and $\tilde K \lo \Theta$ {are} the matrices defined in (\ref{eq:K and K}), then
\begin{align*}
C \breve {\Sigma} \lo {U \Theta}\breve{\Sigma} \lo {U \Theta} \hi * C \hi *=\frac{1}{n \hi 2}\tilde{{K}} \lo {\Theta}{K} \lo U.
\end{align*}
\end{proposition}

Let $\odot$ denote Hadamard product between matrices, and let $1 \lo n$ {be} the $n$-dimensional vector with all entries equal to 1, then we have the following alternative expression for $\|\breve{\Sigma} \lo {U \Theta}\| \loo{HS} \hi 2$:
\begin{equation*}
\|\breve{\Sigma} \lo {U \Theta}\| \loo {HS} \hi 2=\frac{1}{n \hi 2} {1} \lo n \trans ( {K} \lo U\odot\tilde{ {K}} \lo {\Theta}) {1} \lo n.
\end{equation*}
The computation of $ {K} \lo U$ is straightforward. However, for computing $\tilde{ {K}} \lo \Theta$, we  need
\begin{align*}
\tint  \ka \lo \Theta (\hat \Theta \lo i, \theta ) d {P \lo  0} (\theta), \quad  \tint\tint \ka \lo \Theta (\theta, \theta') d ({P \lo 0 \times P \lo 0}) (\theta, \theta').
\end{align*}
We propose to compute these by numerical integration. By Proposition \ref{proposition:n s condition 1}, the $d-1$ components of $\Theta$, namely, $\Theta \lo 1, \ldots, \Theta \lo {d-1}$, are independent with densities
\begin{align*}
f \lo {\Theta \lo j} (\theta \lo j)
=
\begin{cases}
\cos \hi {d-1-j} (\theta \lo j) / \tint \lo {-\pi /2 } \hi {\pi / 2}\cos \hi {d-1-j} (\theta \lo j)  d \theta \lo j,   & \theta \lo j \in (-\pi/2, \pi / 2], \ j = 1, \ldots, d - 2 \\
1/(2 \pi ),  & \theta \lo j \in (-\pi, \pi ], \  j = d - 1.
\end{cases}
\end{align*}
Let
\begin{align*}
\Omega \lo {\Theta \lo j} =
\begin{cases}
(-\pi / 2, \pi / 2 ], & j = 1, \ldots, d - 2 ,\\
( - \pi , \pi ], & j = d - 1.
\end{cases}
\end{align*}
If we choose $\ka \lo \Theta$ to be the product kernel,
\begin{equation}\label{eq:prod-kernel}
\kappa \lo \Theta(\theta,\theta')=\prod \lo {j=1} \hi {d-1}\kappa \lo {\Theta j}(\theta \lo j,\theta \lo j'),
\end{equation}
of which a typical example is the Gaussian radial basis kernel as in (\ref{eq:gaussian kernel}), {then} we have
\begin{align}
\int \lo {\Omega \lo \Theta} \,  \ka \lo \Theta (\hat \Theta \lo i, \theta) d {P \lo 0} (\theta)  = \ali   \prod \lo {j=1} \hi {d-1}  \int \lo {\Omega \lo {\Theta \lo j} } \, \ka \lo {\Theta \lo j} ( \hat \Theta \lo {ij} , \theta \lo j) f \lo {\Theta \lo j} (\theta \lo j) d \theta \lo j, \label{eq:quantities} \\
\int  \lo {\Omega \lo \Theta \times \Omega \lo \Theta} \,  \ka \lo \Theta ( \theta, \theta') d ({P \lo 0 \times  P \lo 0}) (\theta, \theta')
 = \ali  \prod \lo {j=1} \hi {d-2}  \int \lo { \Omega \lo {\Theta \lo j} \times \Omega \lo {\Theta \lo j}} \, \ka \lo {\Theta \lo j} ( \theta \lo j  , \theta \lo j') f \lo {\Theta \lo j} (\theta \lo j)
f \lo {\Theta \lo j} (\theta \lo j') d \theta \lo j d \theta \lo j', \label{eq:quantity}
\end{align}
where, in \eqref{eq:quantities}, $\hat \Theta \lo {ij}$ means the $j$-th element of $\hat \Theta \lo i$.

The quantities in \eqref{eq:quantities} can be computed by the function \texttt{integrate} in R, and the quantity in \eqref{eq:quantity} can be computed by the function \texttt{cubintegrate} in the R package \texttt{cubature} (\cite{cubature}).

\section{Asymptotic distribution}\label{section:null distribution}

In this section we derive the asymptotic distribution of $T \lo n$ under the null hypothesis (\ref{eq:targeted hypothesis}). We use the von-Mises expansion to achieve this purpose; see, for example, \cite{vaart1998asymptotic,fernholz1983vonMC,li2018sufficient}. We first outline the key steps and notations for the von-Mises expansion, tailored for our current application.

\subsection{von-Mises expansion}
Let $\frak{F}$ denote the class of all distributions on $(\Omega \lo X, \ca F \lo X)$. Let $\ca H$ be a generic Hilbert space. Let $T: \frak F \to \ca H$ be a mapping --- such mappings are known as statistical functionals. Endow $\frak F$ with the uniform metric, and $\ca H$ with the metric induced by its inner product. Let $F \lo 0$ be a member of $\frak F$. Then $T$ is Frechet differentiable at $F \lo 0$ if there is a linear operator $A: \frak F \to \ca H$ such that
\begin{align*}
\lim \lo {\| F - F \lo 0 \| \lo {\frak F} \to 0} \frac{ \| T ( F ) - T (F \lo 0) - A (F - F \lo 0) \| \lo {\ca H} }{\| F - F \lo 0 \| \lo {\frak F}} = 0.
\end{align*}
Under the Frechet differentiability, the linear operator $A$ can be calculated using Gateaux derivative: for any $F \in \frak F$, $A(F)$ is simply
\begin{align*}
\lim \lo {\epsilon \to 0} \frac{T((1-\epsilon) F \lo 0 + \epsilon F ) - T (F \lo 0)}{\epsilon}.
\end{align*}
Let $x$ be a member of $\Omega \lo X$, and $\delta \lo x$ be the Dirac measure at $x$. Then the mapping $x \mapsto A(\delta \lo x)$ from $\Omega \lo X$ to $\ca H$ is called influence function of $T$. We write $A ( \delta \lo x)$ as  $T \hi \star (x)$. Note that we use $\star$ to indicate influence function and $*$ to indicate the adjoint operator. Both notations will be used heavily in our exposition. The key result that we use is this: if $X \lo 1, \ldots, X \lo n$ are i.i.d. from $F \lo 0$ and if $T : \frak F \to \ca H$ is Frechet differentiable at $F \lo 0$, then
\begin{align}\label{eq:functional delta}
\sqrt n [ T (F \lo n) - T (F \lo 0)] \cid N(0, \Gamma),
\end{align}
where $\Gamma \in \ca H \otimes \ca H$ is the linear operator {
\begin{align*}
E ( T \hi \star (X) \otimes T \hi \star (X) ) - E ( T \hi \star (X) )  \otimes E ( T \hi \star (X) )=E ( T \hi \star (X) \otimes T \hi \star (X) ) .
\end{align*}}
This fact is known as the {$\delta$-method} for statistical functionals. In our case, $\ca H$ will be the tensor product space $\ca H \lo U \otimes \ca H \lo \Theta$ introduced earlier. From the above discussion we see that the key to computing the asymptotic normal distribution of $\sqrt n ( T (F \lo n) - T (F \lo 0))$ is to compute the influence function $T \hi \star (x)$. Next, we define our statistical functional.

\vspace{-10pt}

\subsection{Statistical functional for testing elliptical distributions}
Now let us consider the statistical functional $T(F)$ in our setting. Since it  is a simple function of $\breve{\Sigma} \lo {U \Theta}$, let us first figure out the statistical functional corresponding to this linear operator.  For any $F \in \frak{F}$, let
\begin{equation*}
\mu(F)=\tint xd F(x),\quad \Sigma(F)=\tint (x-\mu(F))(x-\mu(F))\trans d F(x).
\end{equation*}
Clearly,
$\mu(F \lo 0)=\mu$,     $\mu(F \lo n)=\hat{\mu}$,
$\Sigma(F \lo 0)=\Sigma$,  and    $\Sigma(F \lo n)=\hat{\Sigma}$.
For each $x \in \Omega \lo X$, let
\begin{align*}
W(x,F) =\ali  \Sigma(F) \hi {-1/2}(x-\mu(F)), \\
U(x,F)= \ali \|W(x,F)\|,\\
V(x,F)= \ali W(x,F)/U(x,F),\\
\Theta(x,F)= \ali g(V(x,F)).
\end{align*}
It is easy to see that, when evaluated at $F = F \lo 0$, $U(x,F)$, $V(x,F)$, and $\Theta(x,F)$ reduce to  $  U(x)$, $  V(x)$, and $  \Theta (x)$, and  when evaluated at $F = F \lo n$, they reduce to  $\hat U(x)$, $\hat V(x)$  and $\hat \Theta (x)$.
Our statistical functional of interest is the {Hilbert-Schmidt} norm of the operator
\begin{align}\label{eq:Sigma UT}
\begin{split}
\Sigma \lo {U \Theta}(F)= \ali \tint \kappa \lo U(\cdot,U(x,F))\otimes\kappa \lo {\Theta}(\cdot,\Theta(x,F))d F(x) \\
\ali \quad -\tint \kappa \lo U(\cdot,U(x,F))d F(x)\otimes \tint\kappa \lo {\Theta}(\cdot, \theta )d {P \lo 0}  (\theta).
\end{split}
\end{align}
It is important to note that the last term on the right, $\int\kappa \lo {\Theta}(\cdot, \theta )d {P \lo 0}  (\theta)$,  does not involve the unknown distribution $F$. This is the true expectation determined by the known distribution {$P \lo 0$} in Proposition \ref{proposition:n s condition 1}.
The operator in (\ref{eq:Sigma UT}) can be re-expressed via the centered kernel as
\begin{align}\label{eq:test operator}
\Sigma \lo {U \Theta}(F)= \tint \kappa \lo U(\cdot,U(x,F))\otimes\tilde{\kappa} \lo {\Theta}(\cdot,\Theta(x,F))d F(x).
\end{align}
Clearly, when evaluated at $F=F \lo 0 $, $\Sigma \lo {U \Theta}(F)$ reduces to $\Sigma \lo {U \Theta}$, and when evaluated at $F=F \lo n$, $\Sigma \lo {U \Theta}(F)$ reduces to $\breve{\Sigma} \lo {U \Theta}$.
Our statistical functional of interest is then defined as
\begin{align*}
\frak{F} \to \real, \quad F \mapsto   \| \breve\Sigma \lo {U \Theta} (F) \| \loo{HS} \hi 2.
\end{align*}

\subsection{Derivations of influence functions}

In this subsection we derive the influence functions of statistical functionals involved in $\|\breve \Sigma \lo {U \Theta} (F) \| \loo{HS} \hi 2$. Some of these functionals are of the form $F \mapsto G(x, F)$, which already depends on $x$. To make a distinction with this $x$ and the argument $x$ in $T \hi  \star (x) = A ( \delta \lo x)$, we denote the argument in any influence function by $z$.  Thus, we denote the influence function of the statistical functional $F \mapsto G(x, F)$ as $G \hi \star (x, z)$. That is,
\begin{align*}
G \hi \star (x,z) = [\partial G ( x, (1-\epsilon ) F \lo 0 + \epsilon \delta \lo z) /\partial \epsilon] \lo {\epsilon = 0}.
\end{align*}
We will  refer to the process of deriving $A \hi \star (z)$ from $A (F)$ as the $\star$-operation. The basic rules for the $\star$-operation are given in Proposition 9.2 of \cite{li2018sufficient}. We start with the influence functions about $\mu(F)$, $\Sigma(F)$. The results are given by Lemma 9.1 of \cite{li2018sufficient}, and we reproduce them here for later references.

\begin{lemma}\label{lemma-influence-func-mu-sigma}
If $X$ is integrable, then the influence function of $\mu$ is
\begin{equation*}
\mu \hi \star (z) =z-\mu.\label{eq-mu-star}
\end{equation*}
Furthermore, if $X$ is square integrable, then
\begin{align}
&\Sigma \hi \star (z) =(z-\mu)(z-\mu) \trans  -\Sigma,\nonumber \\
&(\Sigma \hi {-1}) \hi \star (z) =-\Sigma \hi {-1}\Sigma \hi \star (z)\Sigma \hi {-1},\label{eq-sigma-inv-star}\\
&\vec[{(\Sigma \hi {-1/2})} \hi \star (z)]=-(\Sigma \hi {1/2}\otimes\Sigma+\Sigma\otimes\Sigma \hi {1/2}) \hi {-1}\vec[\Sigma \hi {\star} (z)], \label{eq-sigma-inv-sqrt-star}
\end{align}
where, for a matrix $A$ with columns $a \lo 1, \ldots, a \lo m$, $\vec(A)$ denotes the vector $( a \lo 1 \trans, \ldots, a \lo m \trans ) \trans$.
\end{lemma}

We next derive the influence functions for $U(x, F)$, $V(x, F)$, and $\Theta(x, F)$.
For  deriving the influence function of $\Theta(x,F)$, we need the derivative of the polar coordinate transformation, which is given in the next lemma.

\begin{lemma}\label{lemma:dgdv}
Let $v=(v \lo {1},\ldots,v \lo {d} ) \trans  \in\mathbb{S} \hi {d-1}$ and {$\theta=g(v)=(\theta \lo {1},\ldots,\theta \lo {d-1}) \trans \in \Omega \lo \Theta$}. Let $S \lo i$ be {the Euclidean norm of the vector $(v \lo j , \ldots, v \lo d)\trans$ }for $i=1,\ldots,d$.
Then
\begin{equation*}
\frac{\partial g (v)}{ \partial v \trans}=
\begin{pmatrix}
\frac{S \lo 2}{S \lo 1 \hi 2}&
-\frac{v \lo 1v \lo 2}{S \lo 1 \hi 2S \lo 2}&
-\frac{v \lo 1v \lo 3}{S \lo 1 \hi 2S \lo 2}&
\ldots&
-\frac{v \lo 1v \lo {d-2}}{S \lo 1 \hi 2S \lo 2}&
-\frac{v \lo 1v \lo {d-1}}{S \lo 1 \hi 2S \lo 2}&
-\frac{v \lo 1v \lo d}{S \lo 1 \hi 2S \lo 2}\\
0&
\frac{S \lo 3}{S \lo 2 \hi 2}&
-\frac{v \lo 2v \lo 3}{S \lo 2 \hi 2S \lo 3}&
\ldots&
-\frac{v \lo 2v \lo {d-2}}{S \lo 2 \hi 2S \lo 3}&
-\frac{v \lo 2v \lo {d-1}}{S \lo 2 \hi 2S \lo 3}&
-\frac{v \lo 2v \lo d}{S \lo 2 \hi 2S \lo 3}\\
0&
0&
\frac{S \lo 4}{S \lo 3 \hi 2}&
\ldots&
-\frac{v \lo 3v \lo {d-1}}{S \lo 3 \hi 2S \lo 4}&
-\frac{v \lo 3v \lo {d-2}}{S \lo 3 \hi 2S \lo 4}&
-\frac{v \lo 3v \lo d}{S \lo 3 \hi 2S \lo 4}\\
\vdots&
\vdots&
\vdots&
\ddots&
\vdots&
\vdots&
\vdots\\
0&
0&
0&
\ldots&
\frac{S \lo {d-1}}{S \lo {d-2} \hi 2}&
-\frac{v \lo {d-2}v \lo {d-1}}{S \lo {d-2} \hi 2S \lo {d-1}}&
-\frac{v \lo {d-2}v \lo d}{S \lo {d-2} \hi 2S \lo {d-1}}\\
0&
0&
0&
\ldots&
0&
\frac{v \lo d}{S \lo {d-1} \hi 2}&
-\frac{v \lo {d-1}}{S \lo {d-1} \hi 2}
\end{pmatrix}
.
\end{equation*}
\end{lemma}

Based on Lemma \ref{lemma-influence-func-mu-sigma}, we next derive the influence functions for the statistical functionals
\begin{align*}
F \mapsto U(x, F), \quad F \mapsto V(x, F), \quad F \mapsto \Theta(x, F).
\end{align*}

\begin{lemma}\label{lemma:U star V star T star}
Suppose that $F\mapsto U(x, F)$, $F\mapsto V(x, F)$ and $F\mapsto \Theta(x, F)$ are Frechet differentiable at $F_0$. Let
\begin{align}
\begin{split}
A \lo 1 (x)&=-\frac{(x-\mu) \trans \Sigma \hi {-1}}{[(x-\mu) \trans \Sigma \hi {-1} (x-\mu)] \hi {1/2}},\\
A \lo 2 (x)&=-\frac{\left[(x-\mu) \trans \Sigma \hi {-1}\right]\otimes\left[(x-\mu) \trans \Sigma \hi {-1}\right]}{2[(x-\mu) \trans \Sigma \hi {-1}(x-\mu)]\hi {1/2}},\\
B \lo 1 (x)&=-\frac{\Sigma \hi {-1/2}(x-\mu)A \lo 1 (x)}{(x-\mu) \trans \Sigma \hi {-1}(x-\mu)}-\frac{\Sigma \hi {-1/2}}{[(x-\mu) \trans \Sigma \hi {-1}(x-\mu)] \hi {1/2}},\\
B \lo 2 (x)&=-\frac{\Sigma \hi {-1/2}(x-\mu)A \lo 2 (x)}{(x-\mu) \trans \Sigma \hi {-1}(x-\mu)}-\frac{\left[(x-\mu) \trans \otimes I \lo d\right](\Sigma \hi {1/2}\otimes\Sigma+\Sigma\otimes\Sigma \hi {1/2}) \hi {-1}}{[(x-\mu) \trans \Sigma \hi {-1}(x-\mu)] \hi {1/2}}, \\
C \lo i (x)&=[\partial g (V(x))/ \partial v \trans ]B \lo i  (x), \quad i = 1, 2,
\end{split}\label{eq:A B C}
\end{align}
where $\partial g (V(x))/ \partial v \trans$ is the matrix given by Lemma \ref{lemma:dgdv}.
Then, the influence functions of $U(x, F)$, $V(x, F)$ and $\Theta(x, F)$ are
\begin{align}
U \hi \star(x,z)&=A \lo 1 (x)\mu \hi \star(z)+A \lo 2 (x)\vec{[\Sigma \hi \star(z)]}
\label{eq-ustar-final},\\
V \hi \star(x,z)&= B \lo 1 (x)\mu \hi \star(z)+ B \lo 2 (x)\vec{[\Sigma \hi \star(z)]}
\label{eq-vstar-final},\\
\Theta \hi \star(x,z)&=C \lo 1 (x)\mu \hi \star(z)+ C \lo 2 (x)\vec{[\Sigma \hi \star(z)]}.\label{eq-thetastar-final}
\end{align}
\end{lemma}

The next lemma gives the influence function of $\Sigma \lo {U \Theta}(F)$. Henceforth, for a kernel function $\ka (\cdot, t)$, we use $\dot \ka (\cdot, t)$ to denote the partial derivative with respect to the second argument, $\partial \ka (\cdot, t) / \partial t$.

\begin{lemma}\label{lemma:Sigma UT star}
Suppose $F \mapsto \Sigma \lo {U \Theta} (F)$ is Frechet differentiable at $F \lo 0$. Then
\begin{align}\label{eq:Sigma UT star true}
\begin{split}
\Sigma \lo {U \Theta} \hi \star (z)=\ali \kappa \lo U(\cdot,U(z))\otimes\tilde{\kappa} \lo {\Theta}(\cdot,\Theta(z)) - \Sigma \lo {U \Theta} \\
\ali +
E  \{ [\dot \kappa \lo U ( \cdot, U(X))  U \hi \star (X,z)] \otimes\tilde{\kappa} \lo {\Theta}(\cdot,\Theta(X)) \} \\
\ali + E \{ \kappa \lo U(\cdot,U(X))\otimes\dot \kappa \lo \Theta ( \cdot, \Theta(X)) \trans  \Theta  \hi \star (X,z) \},
\end{split}
\end{align}
where $U \hi \star (x,z)$ and $\Theta  \hi \star (x,z)$ are given by Lemma \ref{lemma:U star V star T star}.
\end{lemma}

Even though $\Sigma \lo {U \Theta}$ is always a member of $\ca H \lo U \otimes \ca H \lo \Theta$, it does not follow that $\Sigma \lo {U \Theta} \hi \star (z)$ must also be a member of $\ca H \lo U \otimes \ca H \lo \Theta$. However, to
facilitate computation of (\ref{eq:Sigma UT star true}) at the sample level, we need  $\Sigma \lo {U \Theta} \hi \star (z)$ to be  a member of $\ca H \lo U \otimes \ca H \lo \Theta$. Fortunately, this  can be ensured
if \begin{align}\label{eq:C2 condition}
\ka \lo U \in C \hi 2 ( \Omega \lo U \times \Omega \lo U), \quad \ka \lo \Theta \in C \hi 2 ( \Omega \lo \Theta \times \Omega \lo \Theta),
\end{align}
where, for a set $A \subseteq \real \hi m$,   $C \hi 2(A \times A)$ denotes the set of all real-valued functions on $A\times A$ that are twice differentiable with a bounded Hessian matrix. {The establishment is placed in the Supplementary Materials, which requires a special case of Theorem 1 of  \cite{zhou2008derivative}. }

\subsection{{Asymptotic distribution} of the test statistic}

Based on the influence function of $\Sigma \lo {U \Theta}(F)$ computed in Lemma \ref{lemma:Sigma UT star} and the functional Delta method expressed in (\ref{eq:functional delta}), we can directly write down the asymptotic distribution of $\breve{\Sigma} \lo {U \Theta}$.

\begin{theorem}\label{thm-asymp-normal}
If the statistical functional $F \mapsto {\Sigma} \lo {U \Theta} (F)$ is Frechet differentiable at $F \lo 0$ with respect to the uniform metric in $\frak F$ and conditions in (\ref{eq:C2 condition}) are satisfied, then
\begin{equation*}
\sqrt{n}(\breve{\Sigma} \lo {U \Theta}-\Sigma \lo {U \Theta})\xrightarrow{\mathcal{D}}N(0,\Gamma),
\end{equation*}
where $\Gamma: \ca H \lo U \otimes \ca H \lo \Theta \to \ca H \lo U \otimes \ca H \lo \Theta$ is the operator
\begin{equation*}
\Gamma=E\left[\Sigma \lo {U \Theta} \hi \star(X)\otimes\Sigma \lo {U \Theta} \hi \star(X)\right],
\end{equation*}
and $\Sigma \lo {U \Theta} \hi \star(z)$ is given by (\ref{eq:Sigma UT star true}).
\end{theorem}

Note that the assertion that $\Gamma$ is an operator from $\ca H \lo U \otimes \ca H \lo \Theta$  to $\ca H \lo U \otimes \ca H \lo \Theta$ is a consequence of \eqref{eq:C2 condition}, which guarantees that $\Sigma \lo {U \Theta} \hi \star(z) \in \ca H \lo U \otimes \ca H \lo \Theta$.
Since, under the null hypothesis,  $\Sigma \lo {U \Theta}=0$, we have the following corollary of Theorem \ref{thm-asymp-normal}, which will be  important for sample-level implementation.

\begin{corollary}\label{corollary:under null}
Under the null hypotheses $H \lo 0: \Sigma \lo {U \Theta}=0$ and  the conditions in Theorem \ref{thm-asymp-normal}, we have
\begin{equation}
\sqrt{n}\breve{\Sigma} \lo {U \Theta}\xrightarrow{\mathcal{D}}N(0,\Gamma).\label{eq-asymp-dist-null-hs}
\end{equation}
where $
\Gamma=E\left[\Sigma \lo {U \Theta} \hi \star(X)\otimes\Sigma \lo {U \Theta} \hi \star(X)\right]$ and $\Sigma \lo {U \Theta} \hi \star(z)$ as given in (\ref{eq:Sigma UT star true}) but without the $\Sigma \lo {U \Theta}$ term.
\end{corollary}

Applying continuous mapping theorem to \eqref{eq-asymp-dist-null-hs} by taking the squared Hilbert-Schmidt norm, we have the following corollary.

\begin{corollary}\label{corollary: weighted chisq}
Under the null hypotheses $H \lo 0: \Sigma \lo {U \Theta}=0$ and  the conditions in Theorem \ref{thm-asymp-normal}, we have
\begin{equation}
n\|\breve{\Sigma} \lo {U \Theta}\| \loo{HS}  \hi 2\xrightarrow{\mathcal{D}}\tsum \lo {i=1} \hi \infty\lambda \lo iZ \lo i \hi 2,\label{eq-asymp-dist-null-hs-sq}
\end{equation}
where $Z \lo 1,Z \lo 2,\ldots$ are i.i.d. standard normal random variables, and $\lambda \lo 1,\lambda \lo 2,\ldots$ are eigenvalues of $\Gamma$ in Corollary \ref{corollary:under null}.
\end{corollary}

Note that the null distribution in \eqref{eq-asymp-dist-null-hs-sq} only depends on the eigenvalues of $\Gamma$, which gives us the chance of not having to save the whole $\Gamma$ on sample-level implement.
{
Based on Theorem \ref{thm-asymp-normal}, we can also derive the asymptotic distribution of $\sqrt{n} ( \| \breve \Sigma \lo {U \Theta}  \| \loo {HS} \hi 2 - \|   \Sigma \lo {U \Theta}  \| \loo {HS} \hi 2 )$ under the alternative hypothesis.

\begin{corollary}\label{cor:alternative-distribution}
Suppose that the conditions of Theorem \ref{thm-asymp-normal} hold. Then, under the alternative hypothesis, i.e., $\Sigma \lo {U \Theta} \ne 0$, we have
\begin{align*}
\sqrt{n} ( \| \breve \Sigma \lo {U \Theta}  \| \loo {HS} \hi 2 - \|   \Sigma \lo {U \Theta}  \| \loo {HS} \hi 2 ) \xrightarrow{\mathcal{D}} N (0, 4 E (\langle \Sigma \lo {U \Theta} \hi \star (X), \Sigma \lo {U \Theta} \rangle \loo {HS} \hi 2 )).
\end{align*}
\end{corollary}

Corollary \ref{cor:alternative-distribution} also implies that $\sqrt{n} \| \breve{\Sigma} \lo {U \Theta} \| \loo {HS} \hi 2  \xrightarrow{P} \infty$ under the alternative hypothesis, which verifies the consistency of our test.
}

{
\subsection{Local power analysis}
In this subsection, we derive the local alternative distribution, which can be used to compute the local power. The idea of the proof is similar to  Theorem 13 of \cite{gretton2012kernel}.

\begin{theorem} \label{thm:local power}
Suppose
\begin{enumerate}
\item $\Gamma$ has spectral decomposition $\tsum \lo {j=1} \hi \infty \lambda \lo j (v \lo j \otimes v \lo j)$, where $v \lo 1, v \lo 2, \ldots$ is an orthonormal basis in $\ca H \lo U \otimes \ca H \lo \Theta$;
\item $\Sigma \lo 1$ is a fixed linear operator in $\ca H \lo U \otimes \ca H \lo \Theta$ with expansion $\tsum \lo {j=1} \hi \infty \sigma \lo j  v \lo j$ and $\| \Sigma \lo 1 \| \loo {HS} = c > 0$.
\end{enumerate}
Then, under the local alternative hypothesis
$
H \lo 1 \hiii n: \Sigma \lo {U \Theta} = n \hi {-1/2} \Sigma \lo 1,
$
we have
\begin{align*}
n\|\breve{\Sigma} \lo {U \Theta} \| \loo {HS}  \hi 2  \xrightarrow{\mathcal{D}} \tsum \lo {j=1} \hi {\infty} \lambda \lo j \tilde{Z} \lo j \hi 2
\end{align*}
where $\tilde{Z} \lo j$ are independent $N(\sigma \lo j / \sqrt{\lambda \lo j}, 1) $ random variables.
\end{theorem}

Using this theorem, we  calculate the local power of our test as
\begin{align*}
P \left( n\|\breve{\Sigma} \lo {U \Theta} \| \loo {HS}  \hi 2 > s \right) \to P \left( \tsum \lo {j=1} \hi {\infty} \lambda \lo j \tilde{Z} \lo j \hi 2 > s  \right).
\end{align*}
}
\vspace{-15pt}

\section{Approximating the asymptotic null distribution}\label{section:implement null distribution}

\subsection{Outline of the problem and notations}\label{subsection:outline}
In this section we approximate the asymptotic distribution of $T \lo n$, which is $\tsum \lo {i=1} \hi \infty \lambda \lo i Z \lo i \hi 2$, where $\lambda \lo 1, \lambda \lo 2, \ldots$ are the eigenvalues of $\Gamma$, and $Z \lo 1, Z \lo 2, \ldots$ are i.i.d. $N(0,1)$. The operator $\Gamma$ is estimated by substituting, wherever possible, the expectation $E$ by the sample average $E \lo n$ in the expression $E ( \Sigma \lo {U \Theta} \hi \star \otimes \Sigma \lo {U \Theta} \hi \star )$. Denote  the estimate of $\Gamma$ as $\hat \Gamma$. We use the eigenvalues $\hat \lambda \lo i$ of $\hat \Gamma$ to estimate $\lambda \lo i$ in the asymptotic distribution
 $\tsum \lo {i=1} \hi \infty \lambda \lo i Z \lo i \hi 2$.
We then use the plug-in estimate $\sum \lo {i=1} \hi \infty \hat \lambda \lo i Z \lo i \hi 2$ to approximate the asymptotic distribution of $T \lo n$.

In the following, for an integer $m$, we use $[m]$ to represent the set $\{1, \ldots, m \}$. If  $A$ is a set and $a: [m] \to A$ is a function, we use $a \lo {[m]}$ to represent the vector $(a \lo 1, \ldots, a \lo m ) \trans$. Furthermore, for a function $f$ defined on $A$, and a function $a: [m] \to A$, we use $f ( a \lo {[m]})$ to denote the vector $(f (a \lo 1), \ldots, f ( a \lo m ) ) \trans$. This notation also extends to functions involving other variables. For example, $\ka ( \cdot, a \lo {[m]})$ represents the vector of functions $(\ka (\cdot, a \lo 1), \ldots, \ka (\cdot, a \lo m ) )\trans$; and $f(x,y,z\lo {[m]})$ represents the vector $(f (x,y, z \lo 1), \ldots, f (x,y, z \lo m ) )\trans$.

One of the advantages of the RKHS is that we can use   the representer theorem (see, for example, \cite{scholkopf2001generalized}) to turn an infinite dimensional problem into a finite dimensional problem. Let $\ca H$ be the RKHS generated by a kernel $\kappa: A \times A \to \real$. Suppose  our statistical procedure relies on functions in $\ca H$ only through $f (a \lo 1), \ldots, f (a \lo m)$ for $a \lo 1, \ldots, a \lo m \in A$. Then we can, without loss of generality, restrict our attention to $\hat {\ca H} \subseteq \ca H$, where $\hat {\ca H}$  is spanned by $\ka (\cdot, a \lo 1), \ldots, \ka (\cdot, a \lo m)$. This is because there is always a function $\hat f  \in \hat {\ca H}$ such that $f (a \lo i)= \hat f (a \lo i)$ for $i = 1, \ldots, m$. In fact, let
\begin{align*}
\hat f = t \lo 1 \ka (\cdot, a \lo 1) + \cdots + t \lo m \ka (\cdot, a \lo m) = t \lo {[m]} \trans \ka (\cdot, a \lo {[m]}).
\end{align*}
Up on solving the equation $\hat f ( a \lo {[m]} ) = f ( a \lo {[m]})$, we have $t \lo {[m]} = K \inv f ( a \lo {[m]})$, where $K$ is the Gram matrix $\{\ka (a \lo i, a \lo j) \} \lo {i,j=1} \hi m$. For this reason, in this section, we will reset $\ca H \lo U$ and $\ca H \lo \Theta$ to be the RKHS spanned by
\begin{align*}
\ca B \lo U = \{ \ka \lo U (\cdot, \hat U \lo i): i = 1, \ldots, n \}, \quad \ca B \lo \Theta = \{ \ka \lo \Theta (\cdot, \hat \Theta \lo i): i = 1, \ldots, n \}.
\end{align*}
A basis for $\ca H \lo U \otimes \ca H \lo \Theta$ is
\begin{align*}
\{  \ka \lo U (\cdot, \hat U \lo i) \otimes \ka \lo \Theta (\cdot, \hat \Theta \lo j): \, i,j = 1, \ldots, n \}.
\end{align*}
Let $C \lo {U \Theta}: \ca H \lo U \otimes \ca H \lo \Theta \to \real \hi {n \hi 2}$ be the coordinate mapping that takes $\ka \lo U ( \cdot, \hat U \lo i) \otimes \ka \lo \Theta ( \cdot, \hat \Theta \lo j)$ to $e \lo i \otimes e \lo j$ in $\real \hi {n \hi 2}$. Then, it is easy to see that $C \lo {U \Theta} = C \lo U \otimes C \lo \Theta$, where $C \lo U$ takes $\ka \lo U (\cdot, \hat U \lo i)$ to $e \lo i$ and $C \lo \Theta$ takes $\ka \lo \Theta (\cdot, \hat \Theta \lo j)$ to $e \lo j$.
Let $K \lo {U \Theta}$ be the $\real \hi {n \hi 2 \times n \hi 2}$ Gram matrix whose $((i,j),(i', j'))$th entry is
\begin{align*}
(K \lo {U \Theta} ) \lo {(i,j), (i',  j')} = \ali  \langle  \ka \lo U (\cdot, \hat U \lo i) \otimes \ka \lo \Theta (\cdot, \hat \Theta \lo j), \ka \lo U (\cdot, \hat U \lo i') \otimes \ka \lo \Theta (\cdot, \hat \Theta \lo  j') \rangle \lo {\ca H \lo U \otimes \ca H \lo \Theta} \\
= \ali \ka \lo U (\hat U \lo i, \hat U \lo i') \ka \lo \Theta ( \hat \Theta \lo j, \hat \Theta \lo  j') = ( K \lo U ) \lo {ii'} ( K \lo \Theta ) \lo {jj'}.
\end{align*}
Thus, in matrix notation,
\begin{align*}
K \lo {U \Theta} = K \lo U \otimes K \lo \Theta,
\end{align*}
where $\otimes$ is the Kronecker product between matrices. As discussed earlier, the eigenvalues of $\hat \Gamma$ {are} the same as the eigenvalues of
\begin{align}\label{eq:eigen problem tensor}
K \lo {U \Theta} \hi {1/2} C \lo {U \Theta}  \hat \Gamma C  \lo {U \Theta}  \hi * K  \lo {U \Theta}  \hi {-1/2} =
(K \lo U \otimes K \lo \Theta)  \hi {1/2} (C \lo {U} \otimes C \lo \Theta)  \hat \Gamma ( C  \lo U \otimes C \lo \Theta)   \hi *(K \lo U \otimes K \lo \Theta)   \hi {-1/2}.
\end{align}
So it all boils down to computing the coordinate  $ (C \lo {U} \otimes C \lo \Theta)  \hat \Gamma ( C  \lo U \otimes C \lo \Theta)   \hi *$.

\subsection{Estimation of $\Gamma$}

By Corollary \ref{corollary:under null}, under $U \indep \Theta$,  $\Gamma$ is of the form
\begin{align}\label{eq:Sigma UT star}
\begin{split}
\Sigma \lo {U \Theta} \hi \star (z)= \ali \kappa \lo U(\cdot,U(z))\otimes\tilde{\kappa} \lo {\Theta}(\cdot,\Theta(z))  \\
\ali +
\tint \kappa \lo U \hi \star  (x,z) \otimes\tilde{\kappa} \lo {\Theta}(\cdot,\Theta(x))d F \lo 0 (x) + \tint \kappa \lo U(\cdot,U(x ))\otimes{\kappa} \lo {\Theta} \hi \star (x,z) d F \lo 0 (x) \\
\equiv \ali  \Sigma \lo {U\Theta, 1} \hi \star (z) + \Sigma \lo {U\Theta, 2} \hi \star (z)+ \Sigma \lo {U\Theta,3 } \hi \star (z).
\end{split}
\end{align}

Let $\hat A \lo 1 (x)$ and $\hat A \lo 2 (x)$ be as defined in (\ref{eq:A B C}) with $\mu$ and $\Sigma$ replaced by $\hat \mu$ and $\hat \Sigma$; let $\hat B \lo 1 (x)$ and $\hat B \lo 2 (x)$ be as defined in (\ref{eq:A B C}) with $\mu (x)$, $\Sigma (x)$, $A \lo 1 (x)$, and $A \lo 2(x)$ replaced by $\hat \mu (x)$, $\hat \Sigma (x)$, $\hat A \lo 1 (x)$, and $\hat A \lo 2(x)$; let $\hat C \lo 1 (x)$ and $\hat C \lo 2 (x)$ be as defined in (\ref{eq:A B C}) with $B \lo 1 (x)$ and $B \lo 2 (x)$ replaced by $\hat B \lo 1 (x)$ and $\hat B \lo 2 (x)$. Let $\hat \mu \hi \star (z)$ and $\hat  \Sigma \hi \star (z)$ be as defined in Lemma \ref{lemma-influence-func-mu-sigma} with $\mu$ and $\Sigma$ replaced by $\hat \mu$ and $\hat \Sigma$. Let $\hat U \hi \star(x,z)$ and $\hat \Theta \hi \star (x,z)$ be as defined in (\ref{eq-ustar-final}) and (\ref{eq-thetastar-final})  with $A \lo 1 (x)$, $A \lo 2 (x)$, $C \lo 1 (x)$, $C \lo 2 (x)$, $\mu \hi \star (x)$, and $\Sigma \hi \star (x)$ replaced by $\hat A \lo 1 (x)$, $\hat A \lo 2 (x)$, $\hat C \lo 1 (x)$, $\hat C \lo 2 (x)$, $\hat \mu \hi \star (x)$, and $\hat \Sigma \hi \star (x)$.

To approximate $\Sigma \lo {U \Theta, 1} \hi \star (z)$ in  (\ref{eq:Sigma UT star}), we {replace $U (z)$ and $\Theta (z)$ by $\hat U (z)$ and $\hat \Theta (z)$}, and replace $\tilde \ka \lo \Theta (\cdot, \theta)$ by
\begin{align*}
\tilde \ka \lo \Theta \hiii e (\cdot, \theta) = \ka \lo \Theta(\cdot, \theta) - n \inv \tsum \lo {i=1} \hi n \ka \lo \Theta(\cdot, \hat \Theta \lo i ),
\end{align*}
where the superscript in $\tilde \ka \hiii e$ indicates the word ``empirical'', as we use the sample average { $n\inv \tsum \lo {i=1} \hi n \ka \lo \Theta (\cdot, \hat \Theta \lo i)$ } instead of the population average $E \ka \lo \Theta (\cdot, \Theta) $ to center the kernel $\ka \lo \Theta$. So $\Sigma \lo {U \Theta, 1} \hi \star (z)$ is approximated by
\begin{align*}
\hat \Sigma \lo {U \Theta, 1} \hi \star (z)=\kappa \lo U(\cdot,\hat U(z))\otimes\tilde{\kappa} \hiii e \lo {\Theta}(\cdot,\hat \Theta(z)).
\end{align*}
To approximate $\Sigma \lo {U \Theta, 2} \hi \star (z)$ {and $\Sigma \lo {U \Theta, 3} \hi \star (z)$} in (\ref{eq:Sigma UT star}), we replace $\kappa \lo U \hi \star  (\cdot, x,z)$ and $\kappa \lo \Theta \hi \star  (\cdot, x,z)$ by
\begin{align*}
\hat \kappa \lo U \hi \star  (\cdot, x,z) = \dot \kappa \lo U ( \cdot, \hat U(x))  \hat U \hi \star (x,z), \quad \hat \kappa \lo \Theta \hi \star  (\cdot, x,z) = \dot \kappa \lo \Theta ( \cdot, \hat \Theta (x)) \trans  \hat \Theta  \hi \star (x,z).
\end{align*}
We replace $\int \cdots d F \lo 0 (x)$ by $E \lo n (\cdots)$. So $\Sigma \lo {U \Theta, 2} \hi \star (z)$  is approximated by
\begin{align*}
\hat \Sigma \lo {U \Theta, 2} \hi \star (z)=n \inv \tsum \lo {i=1} \hi n  \hat  \kappa \lo U \hi \star  (\cdot, X \lo i ,z) \otimes\tilde{\kappa} \lo {\Theta} \hiii e (\cdot,\hat \Theta \lo i ).
\end{align*}
Similarly, $\Sigma \lo {U \Theta, 3} \hi \star (z)$ is approximated by
\begin{align*}
\hat \Sigma \lo {U \Theta, 3} \hi \star (z)=n \inv \tsum \lo {i=1} \hi n   \kappa \lo U(\cdot, \hat U \lo i )\otimes \hat{\kappa} \lo {\Theta} \hi \star (\cdot, X \lo i,z).
\end{align*}
Combining the above results, we approximate $\Sigma \lo {U \Theta} \hi \star (z)$ by
\begin{align}\label{eq:everything}
\begin{split}
\hat \Sigma \lo {U \Theta} \hi \star (z) = \ali \hat \Sigma \lo {U \Theta,1} \hi \star (z) + \hat \Sigma \lo {U \Theta,2} \hi \star (z) +\hat \Sigma \lo {U \Theta,3} \hi \star (z) \\
=\ali  \kappa \lo U(\cdot,\hat U(z))\otimes\tilde{\kappa} \hiii e \lo {\Theta}(\cdot,\hat \Theta(z))
+n \inv \tsum \lo {i=1} \hi n  \hat  \kappa \lo U \hi \star  (\cdot, X \lo i ,z) \otimes\tilde{\kappa} \lo {\Theta} \hiii e (\cdot,\hat \Theta \lo i ) \\
\ali +n \inv \tsum \lo {i=1} \hi n   \kappa \lo U(\cdot, \hat U \lo i )\otimes \hat{\kappa} \lo {\Theta} \hi \star (\cdot, X \lo i,z).
\end{split}
\end{align}
To approximate $\Gamma$, we replace $\Sigma \lo {U \Theta} \hi \star (z)$ by $ \hat \Sigma \lo {U \Theta} \hi \star (z)$ and $\int \cdots d F \lo 0(x)$ by $E \lo n (\cdots)$, as follows
\begin{align*}
\hat \Gamma = n \inv \tsum \lo {i=1} \hi n [ \hat \Sigma \lo {U \Theta} \hi \star (X \lo i) \otimes \hat \Sigma \lo {U \Theta} \hi \star (X \lo i) ].
\end{align*}

\subsection{Approximating the  asymptotic distribution of $n \|\breve{\Sigma} \lo {U \Theta}\| \loo {HS}\hi 2$}

Our goal is to find the eigenvalues of $\hat \Gamma$ so as to form the approximation of $\tsum \lambda \lo i Z \lo i \hi 2$. Recall that $\Gamma$  is a self adjoint linear  operator from $\ca H \lo U \otimes \ca H \lo \Theta$ to $\ca H \lo U \otimes \ca H \lo \Theta$, where, as mentioned earlier,  $\ca H \lo U$ and $\ca H \lo \Theta$ can be regarded as  the finite-dimensional RKHS's spanned by $\ca B \lo U$ and $\ca B \lo \Theta$.
As mentioned at the end of Section \ref{subsection:outline}, we need to compute the matrix (\ref{eq:eigen problem tensor}).

We need to review a coordinate mapping rule in a tensor product space.
Suppose $\ca H \lo 1$ and $\ca H \lo 2$ are finite-dimensional Hilbert spaces of dimension $m \lo 1$ and $m \lo 2$, with bases $\ca B \lo 1 = \{ b \lo 1 \hii 1, \ldots, b \lo {m \lo 1} \hii 1 \}$ and $\ca B \lo 2 =\{ b \lo 1 \hii 2, \ldots, b \lo {m \lo 2} \hii 2 \}$, respectively. Let $C \lo 1: \ca H \lo 1 \to \real \hi m$ and $C \lo 2: \ca H \lo 2 \to \real \hi m$ be coordinate mappings with respect to $\ca B \lo 1$ and $\ca B \lo 2$, respectively. Let $G \lo 1$ and $G \lo 2$ be the Gram matrices for $\ca B \lo 1$ and $\ca B \lo 2$. A basis of the tensor product space $\ca H \lo 1  \otimes \ca H \lo 2$ is
\begin{align*}
\ca B \lo {12} = \{ b \lo {ij} = b \lo i \hii 1 \otimes b \lo j \hii 2: i= 1, \ldots, m \lo 1, j = 1, \ldots, m \lo 2 \}.
\end{align*}
For $f \lo 1 \in \ca H \lo 1$ and $f \lo 2 \in \ca H \lo 2$, the tensor product $f \lo 1 \otimes f \lo 2$ can be viewed in two ways: as a member of $\ca H \lo 1 \otimes \ca H \lo 2$ or as an operator from $\ca H \lo 2$ to $\ca H \lo 1$.
As a member of $\ca H \lo 1 \otimes \ca H \lo 2$, the coordinate of $f \lo 1 \otimes f \lo 2$ with respect to $\ca B \lo {12}$ is
\begin{align*}
(C \lo 1 \otimes C \lo 2) (f \lo 1 \otimes f \lo 2) = (C \lo 1 f \lo 1 ) \otimes ( C \lo 2 f \lo 2 ),
\end{align*}
where $\otimes$ on the right is the Kronecker product between matrices. As a linear operator from $\ca H \lo 2$ to $\ca H \lo 2$, the coordinate of $f \lo 1 \otimes f \lo 2$ with respect to the bases $\{\ca B \lo 1, \ca B \lo 2\}$ is
\begin{align*}
C \lo 1 (f \lo 1 \otimes f \lo 2) C \lo 2 \hi *.
\end{align*}
By \cite{li2018nonparametric}, we have
\begin{align}\label{eq:tensor product map rule}
C \lo 1 (f \lo 1 \otimes f \lo 2) C \lo 2 \hi * = (C \lo 1 f \lo 1) (C \lo 2 f \lo 2)\trans G \lo 2.
\end{align}

Applying (\ref{eq:tensor product map rule}) to $\ca H \lo 1 = \ca H \lo 2 =  \ca H \lo U \otimes \ca H \lo \Theta$, we have
\begin{align*}
\ali (C \lo U \otimes C \lo \Theta)  [ \hat \Sigma \lo {U \Theta} \hi \star (X \lo i) \otimes \hat \Sigma \lo {U \Theta} \hi \star (X \lo i) ] (C \lo U \otimes C \lo \Theta) \hi * \\
\ali \hspace{.5in} = [(C \lo U \otimes C \lo \Theta)  \hat \Sigma \lo {U \Theta} \hi \star (X \lo i)] [ (C \lo U \otimes C \lo \Theta)  \hat \Sigma \lo {U \Theta} \hi \star (X \lo i) ] \trans  (K \lo U \otimes K \lo \Theta).
\end{align*}
Therefore, we can reexpress (\ref{eq:eigen problem tensor}) as
\begin{align}\label{eq:target matrix}
\begin{split}
\ali (K \lo U \otimes K \lo \Theta) \hi {1/2}(C \lo U \otimes C \lo \Theta) \hat \Gamma (C \lo U \otimes C \lo \Theta) \hi * (K \lo U \otimes K \lo \Theta) \hi {-1/2} \\
\ali = n \inv \tsum \lo {i=1} \hi n  \, [(K \lo U \otimes K \lo \Theta) \hi {1/2}(C \lo U \otimes C \lo \Theta)  \hat \Sigma \lo {U \Theta} \hi \star (X \lo i)] [(K \lo U \otimes K \lo \Theta) \hi {1/2} (C \lo U \otimes C \lo \Theta)  \hat \Sigma \lo {U \Theta} \hi \star (X \lo i) ] \trans.
\end{split}
\end{align}
{Let
\begin{align*}
\hat m (X \lo i)= \ali n \hi {-1/2} (K \lo U \otimes K \lo \Theta) \hi {1/2}(C \lo U \otimes C \lo \Theta)  \hat \Sigma \lo {U \Theta} \hi \star (X \lo i), \ \ \hat M = \left( \hat m (X \lo 1), \ldots, \hat m (X \lo n) \right) .
\end{align*}   }
Then the matrix in (\ref{eq:target matrix}) can be rewritten as
\begin{align*}
 (K \lo U \otimes K \lo \Theta) \hi {1/2}(C \lo U \otimes C \lo \Theta) \hat \Gamma (C \lo U \otimes C \lo \Theta) \hi * (K \lo U \otimes K \lo \Theta) \hi {-1/2}
=   \hat M \hat M\trans.
\end{align*}

It remains to calculate $\hat m (X \lo i)$, where the key is to calculate  $(C \lo U \otimes C \lo \Theta)  \hat \Sigma \lo {U \Theta} \hi \star (X \lo i)$. By (\ref{eq:everything}), this is
\begin{align*}
(C \lo U \otimes C \lo \Theta) \hat \Sigma \lo {U \Theta, 1} \hi \star (X_i) +(C \lo U \otimes C \lo \Theta) \hat \Sigma \lo {U \Theta, 2} \hi \star (X_i) +(C \lo U \otimes C \lo \Theta) \hat \Sigma \lo {U \Theta, 3} \hi \star (X_i).
\end{align*}
Reading off $\hat \Sigma \lo {U \Theta, 1} \hi \star (z)$ from (\ref{eq:everything}), we have
\begin{align*}
(C \lo U \otimes C \lo \Theta) \hat \Sigma \lo {U \Theta, 1} \hi \star (X \lo i )
= \ali (C \lo U \otimes C \lo \Theta) [\kappa \lo U(\cdot,\hat U \lo i)\otimes\tilde{\kappa} \hiii e \lo {\Theta}(\cdot,\hat  \Theta \lo i )  ] \\
= \ali  [C \lo U \kappa \lo U(\cdot,\hat U \lo i)]\otimes [ C \lo \Theta\tilde{\kappa} \hiii e \lo {\Theta}(\cdot,\hat \Theta \lo i ) ],
\end{align*}
where
\begin{align*}
C \lo U \kappa \lo U(\cdot,\hat U \lo i )= e \lo i, \quad  C \lo \Theta\tilde{\kappa} \hiii e \lo {\Theta}(\cdot,\hat \Theta \lo i ) = e \lo i - 1 \lo n / n.
\end{align*}
Reading off $\hat \Sigma \lo {U \Theta, 2} \hi \star (z)$ from (\ref{eq:everything}), we have
\begin{align}\label{eq:Sigma 2}
\begin{split}
(C \lo U \otimes C \lo \Theta) \hat \Sigma \lo {U \Theta, 2} \hi \star (X \lo i)
=\ali  n \inv \tsum \lo {j=1} \hi n  (C \lo U \otimes C \lo \Theta)[\hat  \kappa \lo U \hi \star  (\cdot, X \lo j ,X \lo i) \otimes\tilde{\kappa} \lo {\Theta} \hiii e (\cdot,\hat \Theta \lo j )] \\
=\ali  n \inv \tsum \lo {j=1} \hi n [C \lo U \hat  \kappa \lo U \hi \star  (\cdot, X \lo j ,X \lo i)] \otimes [C \lo \Theta\tilde{\kappa} \lo {\Theta} \hiii e (\cdot,\hat \Theta \lo j )].
\end{split}
\end{align}
To calculate {$C \lo U \hat  \kappa \lo U \hi \star  (\cdot, X \lo j ,X \lo i)$}, we now give an expression for $C \lo U (f)$ for any $f \in \ca H \lo U$.
\begin{lemma}\label{lemma:general coordinate}
Suppose  $\ca H$ is  a finite-dimensional RKHS generated by a kernel $\ka$, with basis $\ka (\cdot, a \lo a ), \ldots, \ka (\cdot, a \lo m)$,  $f$ is a member of $\ca H$, and $C$ is the coordinate mapping. Then
\begin{align*}
C (f) = K \inv  f ( a \lo {[k]}),
\end{align*}
where $f ( a \lo {[k]})$  is the vector $(f (a \lo 1), \ldots, f ( a \lo k)) \trans$, and $K$ is the Gram matrix $\{ \ka (a \lo i, a \lo j) \} \lo {i,j=1} \hi m$.
\end{lemma}
\proof
 Since $f \in \ca H$, we have  $f = C  (f) \trans \ka ( \cdot, a \lo {[m]} )$ where $\ka ( \cdot, a \lo {[m]} )$ represents the vector of functions $(\ka ( \cdot, a \lo 1 ), \ldots, \ka ( \cdot, a \lo m)) \trans$. Evaluate this equation at $a \lo 1, \ldots, a \lo m$, we have
\begin{align*}
f ( a \lo {[k]}) = K C (f).
\end{align*}
Solving this equation, we have the desired result. \eop

Applying Lemma \ref{lemma:general coordinate} to (\ref{eq:Sigma 2}), we have
\begin{align*}
\ali (C \lo U \otimes C \lo \Theta) \hat \Sigma \lo {U \Theta, 2} \hi \star (X \lo i)
= n \inv \tsum \lo {j=1} \hi n [ K \lo U \inv \dot \kappa \lo U (\hat U \lo {[n]}, \hat U \lo j ) \hat U \hi \star (X \lo j, X \lo i) ] \otimes ( e \lo j - 1 \lo n / n ) \\
\ali  =n \inv \tsum \lo {j=1} \hi n \{ K \lo U \inv \dot \kappa \lo U (\hat U \lo {[n]}, \hat U \lo j ) [\hat A \lo 1 (X \lo j) \hat \mu \hi \star(X \lo i)+ \hat A \lo 2 (X \lo j)\vec{(\hat \Sigma \hi \star(X \lo i))}] \} \otimes ( e \lo j - 1 \lo n / n ).
\end{align*}
Reading off $\hat \Sigma \lo {U \Theta, 3} \hi \star (z)$ from (\ref{eq:everything}), we have
\begin{align}\label{eq:Sigma 3}
\begin{split}
 (C \lo U \otimes C \lo \Theta)\hat \Sigma \lo {U \Theta, 3} \hi \star (X \lo i)
= \ali n \inv \tsum \lo {j=1} \hi n [ C \lo U  \kappa \lo U(\cdot, \hat U \lo j )] \otimes [C \lo \Theta \hat{\kappa} \lo {\Theta} \hi \star (\cdot, X \lo j,X \lo i)]  \\
= \ali n \inv \tsum \lo {j=1} \hi n e \lo j \otimes [C \lo \Theta \hat{\kappa} \lo {\Theta} \hi \star (\cdot, X \lo j,X \lo i)].
\end{split}
\end{align}
Applying Lemma \ref{lemma:general coordinate} again to (\ref{eq:Sigma 3}), we have
\begin{align*}
\ali (C \lo U \otimes C \lo \Theta) \hat \Sigma \lo {U \Theta, 3} \hi \star (X \lo i) = n \inv \tsum \lo {j=1} \hi n e \lo j \otimes [ K \lo \Theta \inv \dot \kappa \lo \Theta (\hat \Theta \lo {[n]}, \hat \Theta \lo j)\hat \Theta \hi \star (X \lo j, X \lo i)]\\
\ali  =n \inv \tsum \lo {j=1} \hi n e \lo j \otimes \{ K \lo \Theta \inv \dot \kappa \lo \Theta (\hat \Theta \lo {[n]}, \hat \Theta \lo j) [\hat C \lo 1 (X \lo j) \hat \mu \hi \star(X \lo i)+ \hat C \lo 2 (X \lo j)\vec{(\hat \Sigma \hi \star(X \lo i))}] \}.
\end{align*}
To summarize, we have
\begin{align*}
\begin{split}
\hat m (X \lo i)  = \ali n \hi {-1/2} (K \lo U \otimes K \lo \Theta)\hi {1/2}  [ e \lo j \otimes (e \lo j - 1 \lo n / n) ] \\
\ali + n \hi {-1/2} (K \lo U \otimes K \lo \Theta)\hi {1/2}   \left\{ n \inv \tsum \lo {j=1} \hi n [ K \lo U \inv \dot \kappa \lo U (\hat U \lo {[n]}, \hat U \lo j ) \hat U \hi \star (X \lo j, X \lo i) ] \otimes ( e \lo j - 1 \lo n / n ) \right\} \\
\ali + n \hi {-1/2} (K \lo U \otimes K \lo \Theta)\hi {1/2}   \left\{
n \inv \tsum \lo {j=1} \hi n e \lo j \otimes [ K \lo \Theta \inv \dot \kappa \lo \Theta (\hat \Theta \lo {[n]}, \hat \Theta \lo j)\hat \Theta \hi \star (X \lo j, X \lo i)] \right\}.
\end{split}
\end{align*}
We need to calculate the eigenvalues $\hat \lambda \lo 1, \ldots, \hat \lambda \lo n$ of the matrix $\hat M \hat M \trans$. Since this is $n \hi 2 \times n \hi 2$, its eigenvalues are expensive to compute if computed directly. However, it is equivalent to calculate the eigenvalues of $\hat M \trans \hat M$, which is an $n \times n$ matrix. We then use $\tsum \lo {i=1} \hi n \hat \lambda \lo i Z \lo i \hi 2$ to approximate the asymptotic null distribution in (\ref{eq-asymp-dist-null-hs-sq}). We  apply the function \texttt{imhof} in the R package \texttt{CompQuadForm} (\cite{CompQuadForm})   to compute the p-value of the distribution of $\tsum \lo {i=1} \hi n \hat \lambda \lo i Z \lo i \hi 2$. If the p-value is smaller than some prespecified significance level $\alpha$, then we reject $H_0$.

\subsection{Validity and consistency}
Since we have plugged the estimated eigenvalues $\hat{\lambda} \lo {1},\ldots,\hat{\lambda} \lo {n}$ into the asymptotic distribution of $n \|\breve{\Sigma} \lo {U \Theta}\| \lo {HS}\hi 2$, we want to show that this step preserves the validity and consistency of the test.

Notice that $\hat{\Gamma}$ is a consistent estimator of $\Gamma$. Suppose that $\sum \lo {j=1} \hi \infty \lambda \lo {j} \hi {1/2} <\infty$, then, by Theorem 1 of \cite{gretton2009fast}, we have
\begin{align*}
\tsum \lo {j=1} \hi {n} \hat{\lambda} \lo {j} Z \lo {j} \hi {2} \xrightarrow{\mathcal{D}} \tsum \lo {j=1} \hi {\infty} \lambda \lo {j} Z \lo {j} \hi {2},
\end{align*}
which indicates the validity of the test. The consistency is straightforward since, under the alternative hypothesis, $\breve{\Sigma} \lo {U \Theta} \xrightarrow{P} \Sigma \lo {U \Theta}$ which is nonzero, indicating that $n \|\breve{\Sigma} \lo {U \Theta}\| \lo {HS}\hi 2 \xrightarrow{P} \infty$ (see, for example, Section 12.4 of \cite{kokoszka2017introduction}).
\vspace{-10pt}

{
\section{Uniform concentration bounds}\label{section:concentration-bound}
In this section, we develop the concentration bounds for $\big | \|\breve{\Sigma} \lo {U \Theta}\| \loo {HS}-\|{\Sigma} \lo {U \Theta}\| \loo {HS}\big |$ in two cases:  we first considered the simple case where $\mu$ and $\Sigma$ are known, and then considered the general case where $\mu$ and $\Sigma$ are unknown. These bounds also allow us to establish the consistency  of our method when the dimension $d$ goes to infinity with the sample size $n$.
\subsection{Case for known  $\mu$ and $\Sigma$}
We first consider the case where $\mu$ and $\Sigma$ are known, which would be true, for example, for testing a spherical distribution, where $\mu = 0$ and $\Sigma = I \lo p$. In this case  $\hat{U}$ and $\hat{\Theta}$ are replaced by $U$ and $\Theta$, and our test statistic reduces to
\begin{align}\label{eq:sigma-utheta-known}
\breve{\Sigma} \lo {U \Theta} = E \lo n [ \kappa \lo U(\cdot, U )\otimes \kappa \lo \Theta(\cdot, \Theta )] - E \lo n [ \kappa \lo U(\cdot, U )] \otimes E[\kappa \lo \Theta(\cdot, \Theta )].
\end{align}
The next theorem gives the concentration bound for $\|\breve{\Sigma} \lo {U \Theta}\| \loo {HS}$, which is similar to   Theorem 7 in \cite{gretton2012kernel}.
\begin{theorem}\label{thm-concentration-known}
Suppose that $\mu$ and $\Sigma$ are known and $\breve{\Sigma} \lo {U \Theta}$ is defined as (\ref{eq:sigma-utheta-known}). Furthermore,  suppose the kernels $\ka \lo U$ and $\ka \lo \Theta$ are bounded:  $0\le \ka \lo U (u,u') \le M \lo U$ and $0 \le \ka \lo \Theta (\theta,\theta') \le M \lo \Theta$ for all $u,u',\theta,\theta'$. Then,
\begin{align}\label{eq:concentration-known}
P\left(\big | \|\breve{\Sigma} \lo {U \Theta}\| \loo {HS}-\|{\Sigma} \lo {U \Theta}\| \loo {HS}\big | \ge t + 4(M \lo U M \lo \Theta / n ) \hi {1/2}\right)\leq \exp\left(-\frac{ t \hi 2n}{10 M \lo UM \lo \Theta}\right).
\end{align}
\end{theorem}

The proof of Theorem \ref{thm-concentration-known} is placed in the Supplementary Materials. Let $u=\frac{n t \hi 2}{10 M \lo U M \lo \Theta}$. Then, \eqref{eq:concentration-known} is equivalent to
\begin{align}\label{eq-concentration-term4}
P\left(\left|\|\breve{\Sigma} \lo {U \Theta}\| \loo {HS}-\|{\Sigma} \lo {U \Theta}\| \loo {HS}\right| \ge (10 M \lo U M \lo \Theta u /n ) \hi {1/2} + 4 (M \lo U M \lo \Theta / n) \hi {1/2} \right)\le e \hi {-u}.
\end{align}

\def\usphere{\mathbb{S} \hi {d-1}}
\def\sumi{\tsum \lo {i=1} \hi n}

\subsection{Case for unknown $\mu$ and $\Sigma$}

Since we would like our  concentration bounds to reflect the behavior with respect to both the sample size $n$ and the dimension $d$, we need to separate out $d$ and $n$ from any constants in our derivations. To make this explicit, we call a constant that doesn't depend on $d$ or $n$ an {\em absolute constant}. We first make the following assumptions.

\begin{assumption}\label{assumption:lambda lambda}
There exist absolute constants $c \lo 1 > 0$ and $c \lo 2 > 0$ such that
\begin{align*}
c \lo 1 \le \lambda \lo \min (\Sigma) \le \lambda \lo \max (\Sigma) \le c \lo 2.
\end{align*}
\end{assumption}
\vspace{-20pt}
\begin{assumption}\label{assumption:uniform sub Gaussian} The random vector $X$ has a uniform sub-Gaussian distribution in the sense that, there is a constant $\sigma \hi 2$ that is independent of $d$ such that, for all $v \in \sphere$,
\begin{align*}
E (e \hi {\lambda \langle v, X - \mu \rangle } ) \le e \hi {\sigma \hi 2 \lambda \hi 2/ 2}.
\end{align*}
\end{assumption}
\vspace{-20pt}
\begin{assumption}\label{assumption:bounded density} The density of $W$ is bounded by $c \hi d$ for some absolute constant $c$.
\end{assumption}
\vspace{-20pt}
\begin{assumption}\label{assumption:bound-lipschitz}
The kernel functions $\ka \lo U$ and $\ka \lo \Theta$ are bounded and Lipschitz continuous:
\begin{align*}
0\le \ka \lo U (u,u') \le M \lo U,\quad \forall u,u', \quad
0 \le \ka \lo \Theta (\theta, \theta') \le M \lo \Theta,\quad \forall \theta, \theta',
\end{align*}
and
\begin{align}
&&\| \ka \lo U (\cdot, u) - \ka \lo U (\cdot, u') \| \le L \lo U |u-u'|,\quad \forall u,u',\label{eq-Lip-u}\\
&&\| \ka \lo \Theta (\cdot,g(v)) - \ka \lo \Theta (\cdot,g(v')) \| \le d\inv L \lo V \| v - v ' \| ,\quad \forall v, v'.\label{eq-Lip-v}
\end{align}
\end{assumption}

\vspace{-10pt}

In the subsequent discussions we will carefully track the indices of the absolute constants as they will eventually appear in the same expression. Let
\begin{align}\label{eq:f1 f2 f3 f4}
\begin{split}
\ali f \lo 1 (n,d,u) = c \lo 1 \sqrt{\frac{d[\log(2d)+u]}{n}}, \\
\ali f \lo 2 (n,d,u) = c \lo 2  d   \max\left\{\left( {\frac{d+u}{n}} \right) \hi {1/4},\left(\frac{d+u}{n}\right) \hi {1/2}\right\}  + c \lo 3 \sqrt d, \\
\ali f \lo 3 (n,d,u) = c \lo 4 d \sqrt{\frac{ d+u }{n}},\\
\ali f \lo 4 (n,u) = c \lo {9} \sqrt{\frac{u}{n}} + c \lo {10} \sqrt{\frac{1}{n}},
\end{split}
\end{align}
where $c \lo 1, c \lo 2, c \lo 3, c \lo 4, c \lo 9, c \lo {10}$ are some positive absolute constants.

\vspace{-10pt}

\begin{theorem}\label{theorem:main tail bound}
Suppose $X$ satisfies Assumption \ref{assumption:uniform sub Gaussian}, and $X \lo 1,\ldots,X \lo n$ are i.i.d samples of $X$. Further suppose $\Sigma$ satisfies Assumption \ref{assumption:lambda lambda}, the density of $W$ satisfies Assumption \ref{assumption:bounded density}, and the kernel functions $\ka \lo U$ and $\ka \lo {\Theta}$ satisfy Assumption \ref{assumption:bound-lipschitz}.  Then, for any $\epsilon>0$, we have
\begin{align*}
\begin{split}
&P \big ( \big |\|\breve{\Sigma} \lo {U \Theta}\| \loo {HS}-\|{\Sigma} \lo {U \Theta}\| \loo {HS}\big | \ge \left[c \lo 7 + 2 c \lo 8 / (\epsilon d) \right] \times\\
&\left[f \lo 3 (n,d,u) f \lo 2 (n,d,u) + c \lo 5 f \lo 1 (n,d,u)\right] + f \lo 4 (n,u) \big ) \le {n  ( c \lo 6 \epsilon ) \hi d } + 7 e \hi {-u},
\end{split}
\end{align*}
where $f \lo 1, f \lo 2, f \lo 3, f \lo 4$ are as defined in \eqref{eq:f1 f2 f3 f4}, and $c \lo 1, \ldots, c \lo {10}$ are some positive absolute constants.
\end{theorem}

The tail bounds in Theorem \ref{theorem:main tail bound} allows us to establish the consistency of the test even when $d$ goes to infinity with $n$, as shown in the next theorem.

\begin{theorem}\label{theorem:high dimensional consistency}
Suppose all conditions in Theorem \ref{theorem:main tail bound} are satisfied. If $\log n \prec d \prec n \hi {1/4}$, then
\begin{align*}
\|\breve{\Sigma} \lo {U \Theta}\| \loo {HS} \cip \|{\Sigma} \lo {U \Theta}\| \loo {HS}  .
\end{align*}
\end{theorem}
}

\vspace{-12pt}

\section{Simulations}\label{section:simulation}

In this section we present some simulation results of using $T \lo n$ to test ellipticity under both the null distribution and the alternative distribution.
{In all the simulations, we use Gaussian radial basis kernels for both $U$ and $\Theta$ as given in (\ref{eq:gaussian kernel}).} To select the tuning parameters, we use the criterion in Section 6.4 of \cite{li2018nonparametric}:
\begin{equation*}
\frac{1}{\sqrt{\gamma \lo U}}={n \choose 2} \hi {-1}\sum \lo {i=1} \hi {n-1}\sum \lo {j=i+1} \hi n |\hat{U} \lo i-\hat{U} \lo j |,\quad \frac{1}{\sqrt{\gamma \lo \Theta}}={n \choose 2} \hi {-1}\sum \lo {i=1} \hi {n-1}\sum \lo {j=i+1} \hi n \|\hat{\Theta} \lo i-\hat{\Theta} \lo j \| .
\end{equation*}

\subsection{Results under the null distribution} \label{sec: simulation null}

In this subsection we perform simulations under the null distribution. We consider scenarios consisting of  different sample sizes and dimensions:
\begin{align}\label{eq:n and d}
(n, d)  \in \{500, 1000\} \times \{3,4,5,6,10,15,20\}.
\end{align}
For each $(n,d)$, we generate $T=100$ datasets as follows. We first generate the mean vector $\mu$ from $N(0,\sigma \lo \mu \hi 2 I \lo d)$, where $\sigma \lo \mu \hi 2=100$, and the covariance matrix $\Sigma$ using R function \texttt{genPositiveDefMat} in the package \texttt{clusterGeneration} (\cite{clusterGeneration}) under its default settings. We then simulate $X \lo 1,\ldots,X \lo n$ as i.i.d. samples from $N(\mu,\Sigma)$.

\begin{figure}[htp]
\centering
\includegraphics[width=\textwidth]{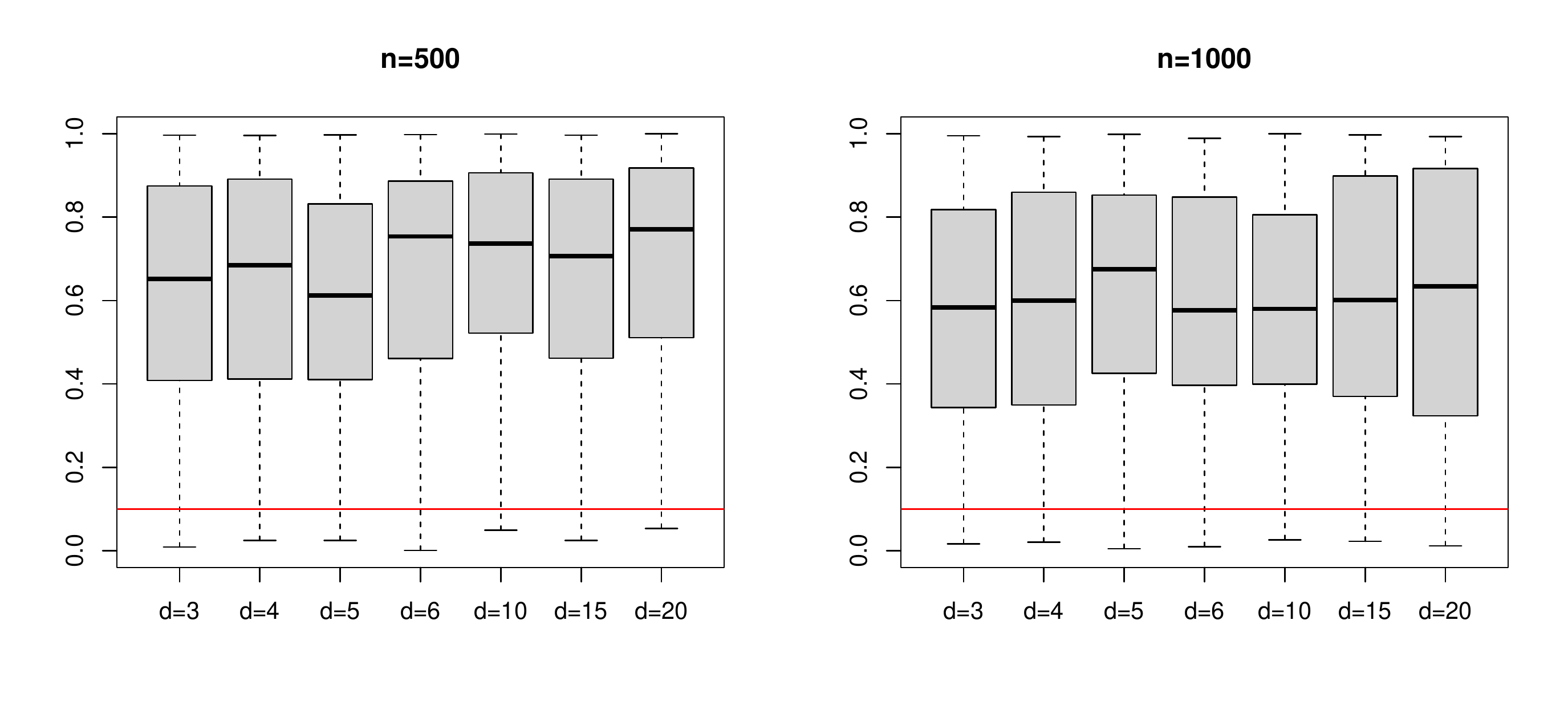}\vspace{-20pt}
\caption{Boxplots of p-values under the null hypothesis. The red line represents 0.1.}
\label{fig-simu-null-pval}
\end{figure}

Based on the samples $X \lo 1, \ldots, X \lo n$ we compute the test statistic $ n \| \breve{\Sigma} \lo {U \Theta} \| \hi 2 \loo {HS}$ and compute the p-value based on our test. One side-note is that we add a small number, $\epsilon=10 \hi {-6}$, to the diagonal of a matrix whenever we compute its inverse, square-root or eigenvalues. This is done to avoid numerical instability that can happen in  some extreme samples.
Figure \ref{fig-simu-null-pval} shows the boxplots of p-values under different combinations of $(n,d)$ in (\ref{eq:n and d}), and {the empirical type-I errors among 100 experiments at $\alpha=0.1$ are summarized in Table \ref{tab:results-null-hypothesis}.}

\vspace{-20pt}

\vspace{-20pt}
\begin{table}[htbp]
\centering
\caption{Empirical type-I errors (in percentage) at $\alpha=0.1$.}\label{tab:results-null-hypothesis}
\begin{tabular}{rrrrrrrr}
\hline
$n$  & $d=3$ & $d=4$ & $d=5$ & $d=6$ & $d=10$ & $d=15$ & $d=20$\\
\hline
500 & 6 & 4 & 5 & 4 & 5 & 1 & 2\\
1000 & 7 & 6 & 3 & 8 & 2 & 7 & 7\\
\hline
\end{tabular}
\end{table}

As we can see {from Figure \ref{fig-simu-null-pval},} the bulk of the p-values are quite large for all combinations of $n$ and $d$,  indicating the elliptical distribution hypothesis is not rejected in the great majority of cases. {In Table \ref{tab:results-null-hypothesis}, the empirical type-I errors are all below the significance level, which is consistent with the validity of our test.}

\subsection{Results under alternative distributions} \label{sec: simulation alternative}

We consider alternative distributions with different degrees of departure from the elliptical distribution. We first generate $Z \lo k$ independently from $N(0, 4)$, and set $Z = (Z \lo 1, \ldots, Z \lo d) \trans$. We then randomly select a subset $J$ of $\{1, \ldots, d \}$ of cardinality $\lceil d/3 \rceil$ and, for each $k \in J$, we replace $Z \lo k$ by $W \lo k  - df$, with $W \lo k $ generated from $\chi \hi 2 (df)$. We denote the resulting random vector by $\tilde Z$.
We then construct $X$ by
\vspace{-10pt}
\begin{align*}
X=\mu+\Sigma \hi {1/2}\tilde{Z},
\end{align*}
where $\Sigma$ is generated in the same way as it was in the null distribution case, and $\mu = (\mu \lo 1, \ldots, \mu \lo d)\trans$ is generated from
\vspace{-10pt}
\begin{align*}
\mu \lo k=40\beta \lo k-20, \quad k = 1, \ldots, d,
\end{align*}
where each $\beta \lo k$ is independently generated from $\mathrm{Beta}(0.5,0.5)$. That is, each $\mu \lo k$ is a rescaled and centered Beta variable. We take the degrees of freedom $df$ to be $2$ or $4$, with $df = 2$ representing stronger departure from ellipticity.

\begin{figure}[htp]
\centering
\includegraphics[width=\textwidth]{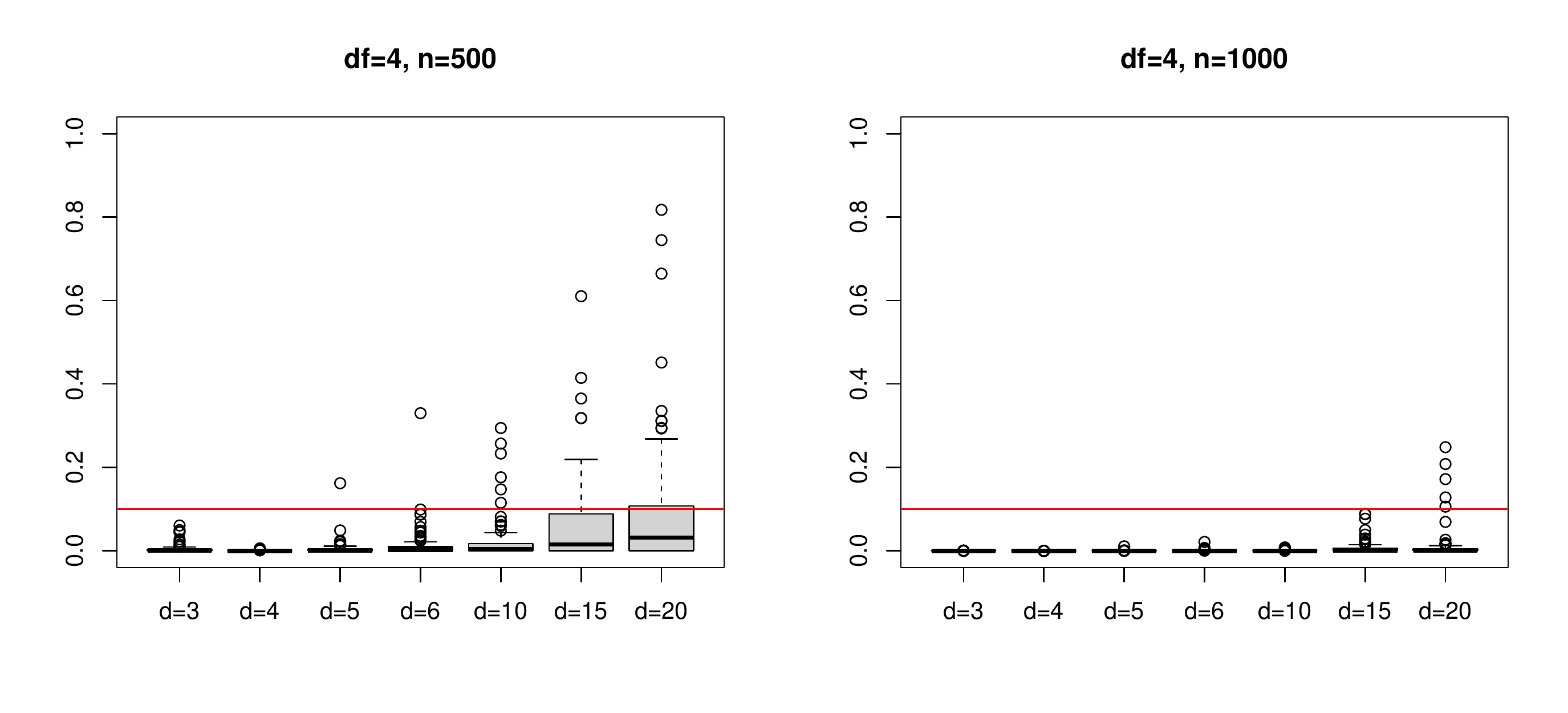}\vspace{-20pt}
\caption{Boxplots of p-values under the alternative hypothesis with $df=4$. The red line represents 0.1.}
\label{fig-simu-alternative-df4-pval}
\end{figure}

\begin{figure}[htp]
\centering
\includegraphics[width=\textwidth]{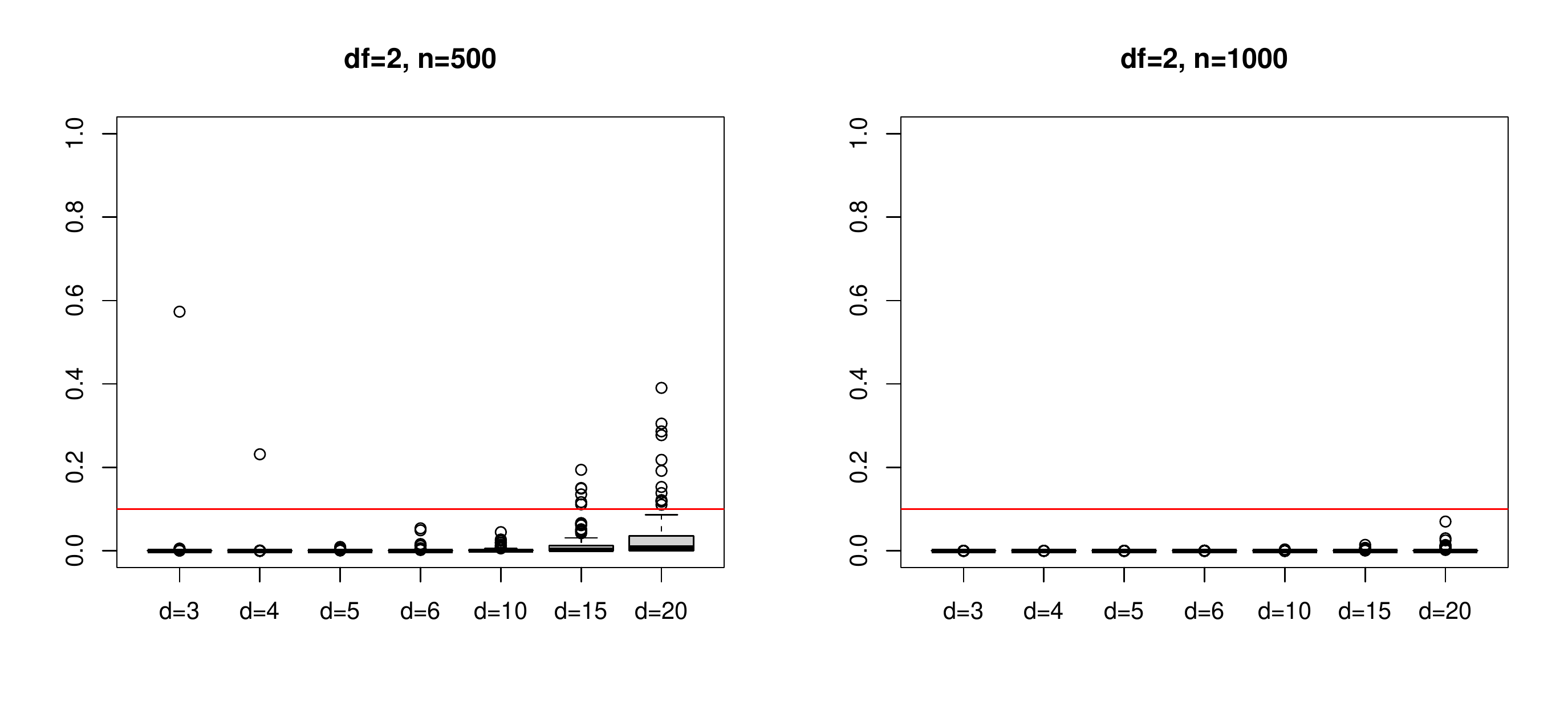}\vspace{-20pt}
\caption{Boxplots of p-values under the alternative hypothesis with $df=2$. The red line represents 0.1.}
\label{fig-simu-alternative-df2-pval}
\end{figure}

We then perform our proposed test on the i.i.d. samples  $X \lo 1, \ldots, X \lo n$ from $X$.  The boxplots in  Figure \ref{fig-simu-alternative-df4-pval} show p-values for $df=4$ with $(n,d)$ in the range (\ref{eq:n and d}); those in  Figure \ref{fig-simu-alternative-df2-pval} show p-values for $df=2$ with $(n,d)$ in the same range. {The numerical values of the empirical powers among 100 experiments at $\alpha=0.1$ are summarized in Table \ref{tab:results-alternative-hypothesis}.}

\vspace{-20pt}

\newpage
\begin{table}[!htbp]
\centering
\caption{Results under the alternative hypothesis}\label{tab:results-alternative-hypothesis}
\begin{tabular}{rrrrrrrrr}
\hline
$n$ & $df$ & $d=3$ & $d=4$ & $d=5$ & $d=6$ & $d=10$ & $d=15$ & $d=20$\\
\hline
500 & 4 & 100 & 100 & 99 & 99 & 94 & 78 & 73\\
500 & 2 & 99 & 99 & 100 & 100 & 100 & 94 & 88\\
1000 & 4 & 100 & 100 & 100 & 100 & 100 & 100 & 95\\
1000 & 2 & 100 & 100 & 100 & 100 & 100 & 100 & 100\\
\hline
\end{tabular}
\end{table}

According to Figure \ref{fig-simu-alternative-df4-pval} and Figure \ref{fig-simu-alternative-df2-pval}, the bulk of the p-values are smaller than 0.1, indicating strong evidence that $X$ does not have an elliptical distribution. This{, together with the empirical powers in Table \ref{tab:results-alternative-hypothesis}, }implies the consistency of our test.
More specifically, in each figure, when the sample size $n$ is fixed, the boxplots becomes slightly higher as $d$ increases, {and the empirical power is lower when $d$ is high}, indicating it is more difficult to detect non-ellipticity in higher dimensions.   This is reasonable because the skewness might be masked  by higher dimensions. If we fix the dimension $d$ and the skewness as represented by $df$, an increase of sample size $n$ from 500 to 1000 makes the p-values more concentrated around 0. Furthermore, comparing Figure \ref{fig-simu-alternative-df2-pval} with Figure \ref{fig-simu-alternative-df4-pval}, we see that the p-values for $df=2$ are smaller than those for $df=4$, indicating that an increase of skewness makes non-ellipticity more detectable.

\section{Application}\label{section:application}

In this section we apply our test to a dataset concerning a Study on the Efficacy of Nosocomial Infection Control (SENIC Project), which is used in \cite{haley1980senic}. We download the dataset from \url{https://users.stat.ufl.edu/~rrandles/sta4210/Rclassnotes/data/textdatasets/KutnerData/Appendix\%20C\%20Data\%20Sets/APPENC01.txt}.
This is one of the datasets from the book \cite{kutner2005applied}.

According to \cite{kutner2005applied}, the data of 113 hospitals during the 1975-76 study period are sampled from the original 338 hospitals surveyed. For a single hospital, there are 11 variables:
\begin{quote}
length of stay, age, infection risk, routine culturing, routine chest x-ray, number of beds, medical school affiliation, region,  average daily census,  number of nurses,  available facilities.
\end{quote}
For detailed information about these variables, see \cite{kutner2005applied}.

We treat the variable ``length of stay'' as the response variable, and check whether other predictor variables are elliptically distributed. We remove the two categorical variables, ``medical school affiliation'' and ``region''. We carry out  our test on the dataset where $n=113$ and $d=8$ at a significance level $\alpha=0.05$, and obtained the  p-value 0.0014, which leads to  rejection of the  elliptical distribution hypothesis.
The original data is as {shown in }Figure 1 in Appendix E of the online Supplementary Material (\cite{supplementary}).

We then perform the Box-Cox transformation on the dataset. Using the method presented in Chapter 7 of \cite{li2018sufficient}, we find the optimal $\lambda$ for the Box-Cox transformation as
\begin{equation*}
1.158,  0.737, -0.105,  0.947, -0.316, -0.105, -0.316,  0.947.
\end{equation*}
The transformed variables are shown in Figure 2 in Appendix E of the online Supplementary Material (\cite{supplementary}). We carry out  our test on the transformed data, and obtain the p-value 0.2454. This is much larger than 0.05, leading us to accept the elliptical distribution hypothesis after the Box-Cox transformation.

{
\section{Discussion}\label{section:discussion}
For our method, the choice of kernels are indeed important. In general, the kernel functions $\kappa \lo U$ and $\kappa \lo \Theta$ should satisfy the following three conditions.

\begin{enumerate}
\item Both $\kappa \lo U$ and $\kappa \lo \Theta$ should be characteristic kernels. This is to guarantee that $\Sigma \lo {U \Theta} = 0$ if and only if $X$ follows an elliptical distribution. If either of them fails to be characteristic, then there will exist cases when $X$ does not follow an elliptical distribution and yet $\Sigma \lo {U \Theta} = 0$ is satisfied, which leads to non-consistency of the test.
\item Both $\kappa \lo U$ and $\kappa \lo \Theta$ should be $C \hi 2$-smooth kernel functions, i.e., the conditions in \eqref{eq:C2 condition} should be satisfied. This is to guarantee that $\Sigma \lo {U \Theta} \hi \star (z)$ is a member of $\mathcal{H} \lo U \otimes \mathcal{H} \lo \Theta$. Also, we need the derivatives $\dot{\kappa} \lo U ( \hat{U} \lo i, \hat{U} \lo j)$ and $\dot{\kappa} \lo \Theta ( \hat{\Theta} \lo i, \hat{\Theta} \lo j)$ in the implementation, which requires differentiability in $\kappa \lo U$ and $\kappa \lo \Theta$.
\item  From a computational perspective, a product-type kernel for $\kappa \lo \Theta$ is preferred. This is because, otherwise, we will need to compute $n$ numerical integrals of dimension $(d-1)$ in equation \eqref{eq:quantities}  and one numerical integral of dimension $(2d-2)$ in equation \eqref{eq:quantity}. The product-type kernels allow us to replace these high-dimension numerical integration  by 1- or 2-dimensional integration.
\end{enumerate}

Clearly, the Gaussian kernel satisfies all three conditions above, so it is preferred. It also works well in our simulation studies.  However, many other kernels also satisfy the three conditions. We can construct a broad class of product kernels that are computationally feasible for our method.

Specifically, for any  $d-1$ kernel functions on $\real \times \real$, say $\kappa \lo {\Theta \lo 1},\ldots,\kappa \lo {\Theta \lo {d-1}}$, their product given by \eqref{eq:prod-kernel} is still a reproducing kernel.  Furthermore,
by Theorem \ref{thm:joint characteristic},    the product of characteristic kernels is still characteristic. For example, we can use the product-type inverse-quadratic (PIQ) kernels defined as follows:
\begin{align}\label{eq:PIQ-kernel}
\kappa \lo U (u,u')= \frac{1}{ 1 + \gamma \lo U (u-u') \hi 2 },\quad
\kappa \lo \Theta (\theta,\theta')= \prod \lo {j=1} \hi {d-1}  \frac{1}{ 1 +\gamma \lo \Theta (\theta \lo j  - \theta \lo j') \hi 2 },
\end{align}
which also satisfy the above three requirements.
In principle, all such kernels can be used for our purpose. This gives us broad choices for kernels.Further discussions on the choice of kernels are given in the Supplementary Materials.
}

%
%

\begin{acks}[Acknowledgments]
The authors would like to thank two referees and an Associate Editor for their insightful comments and suggestions, which helped us greatly in improving this work.
\end{acks}
\begin{funding}
The research of Bing Li is supported in part by the U.S. National Science Foundation (NSF) Grant DMS-2210775 and the U.S. National Institutes of Health (NIH) grant 1 R01 GM152812-01.

%
\end{funding}

\begin{supplement}
\stitle{Supplementary Material for ``A nonparametric test for elliptical distribution based on kernel embedding of probabilities''}
\sdescription{Most of the  proofs, additional simulation studies,  further discussions, and the scatter plot matrices for Section \ref{section:application} can be found in the Supplementary Material.}
\end{supplement}


\bibliographystyle{imsart-nameyear} 
\bibliography{elliptical}       


\end{document}